\documentclass[10pt,reqno]{amsart}
\usepackage[english]{babel}
\usepackage[utf8]{inputenc}
\usepackage{amssymb,amsmath,amsfonts,amsthm,amscd,latexsym,graphicx,epsfig,colordvi}

\usepackage{geometry} 
\let\oldsection\section
\renewcommand{\section}{
  \renewcommand{\theequation}{\thesection.\arabic{equation}}
  \oldsection}
\let\oldsubsection\subsection
\renewcommand{\subsection}{
  \renewcommand{\theequation}{\thesubsection.\arabic{equation}}
  \oldsubsection}

\begin{document}

\title[ROTA-BAXTER OPERATORS OF WEIGHT ZERO ON CAYLEY-DICKSON ALGEBRA]{ROTA-BAXTER OPERATORS OF WEIGHT ZERO ON CAYLEY-DICKSON ALGEBRA}
\author{A.S.Panasenko}
\address{Alexander Sergeevich Panasenko
\newline\hphantom{iii} Novosibirsk State University,
\newline\hphantom{iii} Universitetskiy pr., 1,
\newline\hphantom{iii} 630090, Novosibirsk, Russia
}
\email{a.panasenko@g.nsu.ru}%

\markboth{A.S. Panasenko}
{Rota-Baxter Operators of Weight Zero on Cayley-Dickson algebra}



\maketitle{\small
\begin{quote}
\noindent{\sc Abstract. } All Rota-Baxter operators of weight zero on split octonion algebra over a~field of characteristic not 2 are classified up to conjugation by automorphisms and antiautomorphisms. Thus, the classification of Rota-Baxter operators on composition algebras is finished. There are two descriptions: a~common description over arbitratry field of characteristic not 2 and more accurate description over a~quadratically closed field of characteristic not 2. \medskip

 \end{quote}
}

\section{Introduction}

The Rota-Baxter operator is a~formal generalization of the integration by parts formula \cite{Baxter}. In the work \cite{Semenov} Rota-Baxter operators of nonzero weight appeared independently as solutions of the modified Yang-Baxter equation. At present, applications of Rota-Baxter operators to various areas of algebra are known \cite{Aguiar, Double, Double2}. We are interested in the problem of classifying Rota-Baxter operators on various algebras, especially simple finite-dimensional ones. Descriptions of Rota-Baxter operators are known on $M_2(F)$ over an~algebraically closed field \cite{BGP}, the simple Lie algebra $\mathrm{sl_2(\mathbb{C})}$ \cite{Pan,Pei}, the simple Jordan superalgebra $D_t$ over an~algebraically closed field of characteristic 0 \cite{Bolotina}, $K_3$~\cite{BGP}.

Descriptions of Rota-Baxter operators of nonzero weight are known on a~simple Jordan algebra of a~bilinear form of odd dimension, $M_2(F)$, $K_3$ (all in \cite{BGP}), $M_3(\mathbb{C})$ \cite{M1,M2,M3}.

Composition algebras arose within the framework of a~generalized formulation of the Hurwitz problem \cite{Zhevl}. They can be of two types: split and division algebras. In the paper \cite{BGP} it is proved that all Rota-Baxter operators on a~division quadratic algebra are trivial. Each composition algebra is quadratic. Over a~fixed field $F$ of characteristic not 2 there exist only three split composition algebras: the direct sum of two fields $F$, the matrix algebra $M_2(F)$, and the split octonions over F. In the article \cite{An} Rota–Baxter operators on the direct sum of two fields were described  (later, the description was generalized to the direct sum of a~finite number of fields, \cite{GubarevField}). In \cite{BGP}, Rota–Baxter operators on the algebra of second-order matrices over an~algebraically closed field were described. Thus, the problem of classification of Rota–Baxter operators on composition algebras was reduced to the problem of describing Rota–Baxter operators on split octonions.

Octonions are one of the most famous and well-studied non-associative algebraic systems. Octonions have many mathematical (\cite{Springer}) and physical (\cite{Dixon}) applications. A~brief introduction to octonions can be found in (\cite{Zhevl}, Chapter~2).

In the recent paper \cite{Octonions}, all subalgebras of octonions were described up to automorphism. This description forms the basis for our work. All classifications in this paper are made up to automorphisms and antiautomorphisms of split octonions.

In this paper we complete the classification of Rota-Baxter operators of weight zero on split octonions over a~field of characteristic not 2. The paper is organized as follows. In Section 2 we define several automorphisms and antiautomorphisms on split octonions. In addition, we refine the description of RB-operators on the algebra of second-order matrices from \cite{BGP}, extending the restriction to any field of characteristic not 2. In Section 3 we describe RB-operators of weight zero on split octonions with one-dimensional image. In Section 4 we describe RB-operators of weight zero on split octonions with two-dimensional image. In Section 5 we describe RB-operators of weight zero on split octonions with three-dimensional image. In Section 6 we describe RB-operators of weight zero on split octonions with four-dimensional image. In Section 7 we collect the auxiliary results of the previous sections into a~final result. We formulate two versions of the result: for an~arbitrary field of characteristic not 2 and for a~quadratically closed field of characteristic not~2.

\section{Preliminaries}

In this article we fix a~field $F$ with characteristic $\neq 2$.

An antiautomorphism $\varphi: A \to A$ of an~algebra $A$ is called an~\textbf{\textit{involution}} if $\varphi^2 = \mathrm{id}|_{A}$.

\medskip Let $\mathbb{O}=M_2(F) + vM_2(F)$ be the split Cayley-Dickson algebra. It has the following multiplication table:
\[a\cdot b = ab, \quad a\cdot vb = v(\overline{a}b), \quad va\cdot b = v(ba), \quad va\cdot vb = b\overline{a},\]
where $x\cdot y$ is a~multiplication in $\mathbb{O}$ for any $x,y\in \mathbb{O}$, $ab$ is a~multiplication in $M_2(F)$ for any $a,b\in M_2(F)$, $\overline{a}$ is a~symplectic involution in $M_2(F)$, i.e.
\[\overline{\begin{pmatrix}
    a_{11} & a_{12} \\
    a_{21} & a_{22}
\end{pmatrix}} = 
\begin{pmatrix}
    \phantom{-}a_{22} & -a_{12} \\
    - a_{21} & \phantom{-}a_{11}
\end{pmatrix}.\]
An involution $x\to \overline{x}$ in $\mathbb{O}$ is called \textbf{\textit{classical involution}}.

An algebra $\mathbb{O}$ has the following basis
\[e_{11},e_{12},e_{21},e_{22},ve_{11},ve_{12},ve_{21},ve_{22}.\]

If $\varphi:M_2(F)\to M_2(F)$ is (anti)automorphism then it can be extended to (anti)automorphism $\overline{\varphi}:\mathbb{O}\to\mathbb{O}$ by (\cite{Octonions}, Lemma 4.1). We will use this fact without further mentions.

If $B$ is a~subalgebra in $\mathbb{O}$ then we will call $B$ \textit{\textbf{a unital}} subalgebra if $B$ contains a~unit of $\mathbb{O}$.

The proofs of the following propositions are straightforward.

\medskip\textbf{Proposition 1.} \textit{Let $\varphi:\mathbb{O}\to \mathbb{O}$ be a~linear map such that $\varphi^2=\mathrm{id}$, $\varphi(ve_{12})=-ve_{21}$, $\varphi(ve_{11})=ve_{22}$, $\varphi(e_{11})=e_{11}$, $\varphi(e_{12})=e_{21}$, $\varphi(e_{22})=e_{22}$. Then $\varphi$ is an~involution on the algebra $\mathbb{O}$.}

\medskip\textbf{Proposition 2.} \textit{Let $\varphi:\mathbb{O}\to \mathbb{O}$ be a~linear map such that $\varphi(e_{ij})=e_{ij}$ for any $i,j\in\{1,2\}$, $\varphi(ve_{12})=ve_{12}$, $\varphi(ve_{22})=ve_{22}$, $\varphi(ve_{11})=ve_{11}+\alpha ve_{12}$, $\varphi(ve_{21})=ve_{21}+\alpha ve_{22}$ for some $\alpha\in F$. Then $\varphi$ is an~automorphism on the algebra $\mathbb{O}$.}

\medskip\textbf{Proposition 3.} \textit{Let $\varphi:\mathbb{O}\to \mathbb{O}$ be a~linear map such that $\varphi(e_{ij})=e_{ij}$ for any $i,j\in\{1,2\}$, $\varphi(ve_{11})=ve_{11}$, $\varphi(ve_{21})=ve_{21}$, $\varphi(ve_{12})=ve_{12}+\alpha ve_{11}$, $\varphi(ve_{22})=ve_{22}+\alpha ve_{21}$ for some $\alpha\in F$. Then $\varphi$ is an~automorphism on the algebra $\mathbb{O}$.}

\medskip\textbf{Proposition 4.} \textit{Let $\varphi:\mathbb{O}\to \mathbb{O}$ be a~linear map such that $\varphi(e_{ij})=e_{ij}$ for any $i,j\in\{1,2\}$, $\varphi(ve_{12})=-ve_{11}$, $\varphi(ve_{11})=ve_{12}$, $\varphi(ve_{21})=ve_{22}$, $\varphi(ve_{22})=-ve_{21}$. Then $\varphi$ is an~automorphism on the algebra $\mathbb{O}$.}

\medskip\textbf{Proposition 5.} \textit{Let $\varphi:\mathbb{O}\to \mathbb{O}$ be a~linear map such that $\varphi(e_{11})=e_{22}$, $\varphi(e_{12})=ve_{22}$, $\varphi(e_{21})=ve_{11}$, $\varphi(e_{22})=e_{11}$ and $\varphi^2 =\mathrm{id}$. Then $\varphi$ is an~involution on the algebra $\mathbb{O}$.}

\medskip\textbf{Proposition 6.} \textit{Let $\varphi:\mathbb{O}\to \mathbb{O}$ be a~linear map such that $\varphi(e_{ii})=e_{ii}$, $\varphi(ve_{ii})=ve_{ii}$ for any $i\in\{1,2\}$, $\varphi(e_{12}) = \alpha e_{12}$, $\varphi(ve_{12}) = \alpha ve_{12}$, $\varphi(e_{21})=\frac{1}{\alpha} e_{21}$, $\varphi(ve_{21}) = \frac{1}{\alpha} ve_{21}$ for some $0\neq \alpha\in F$. Then $\varphi$ is an~automorphism on the algebra $\mathbb{O}$.}

\medskip\textbf{Proposition 7.} \textit{Let $\varphi:\mathbb{O}\to \mathbb{O}$ be a~linear map such that $\varphi(e_{ii})=e_{ii}$, for any $i\in\{1,2\}$, $\varphi(e_{12}) = \alpha e_{12}$, $\varphi(ve_{11}) = \alpha ve_{11}$, $\varphi(e_{21})=\frac{1}{\alpha} e_{21}$, $\varphi(ve_{22}) = \frac{1}{\alpha} ve_{22}$, $\varphi(ve_{12})=ve_{12}$, $\varphi(ve_{21})=ve_{21}$ for some $0\neq \alpha\in F$. Then $\varphi$ is an~automorphism on the algebra $\mathbb{O}$.}

\medskip\textbf{Proposition 8.} \textit{Let $\varphi:\mathbb{O}\to \mathbb{O}$ be a~linear map such that $\varphi(e_{ij})=e_{ij}$, for any $i,j\in\{1,2\}$, $\varphi(ve_{11}) = \alpha ve_{11}$, $\varphi(ve_{21})=\alpha ve_{21}$, $\varphi(ve_{22}) = \frac{1}{\alpha} ve_{22}$, $\varphi(ve_{12})= \frac{1}{\alpha} ve_{12}$ for some $0\neq \alpha\in F$. Then $\varphi$ is an~automorphism on the algebra $\mathbb{O}$.}

\medskip\textbf{Proposition 9.} \textit{Let $\varphi:\mathbb{O}\to \mathbb{O}$ be a~linear map such that $\varphi(e_{11}) = e_{11} + \alpha ve_{22}$, $\varphi(e_{12}) = e_{12} + \alpha ve_{12}$, $\varphi(e_{21}) = e_{21}$, $\varphi(e_{22}) = e_{22} - \alpha ve_{22}$, $\varphi(ve_{11}) = ve_{11}-\alpha e_{11}+\alpha e_{22}-\alpha^2 ve_{22}$, $\varphi(ve_{12}) = ve_{12}$, $\varphi(ve_{21}) = ve_{21} + \alpha e_{21}$, $\varphi(ve_{22}) = ve_{22}$ for some $\alpha\in F$. Then $\varphi$ is an~automorphism on the algebra $\mathbb{O}$.
}

\medskip\textbf{Proposition 10} \textit{Let $\varphi:\mathbb{O}\to \mathbb{O}$ be a~linear map such that $\varphi(e_{ii}) = e_{ii}$, $\varphi(ve_{ii}) = ve_{ii}$ for $i\in\{1,2\}$, $\varphi(e_{12}) = e_{12}$, $\varphi(ve_{12}) = ve_{12}$, $\varphi(e_{21}) = e_{21}+\alpha ve_{12}$, $\varphi(ve_{21}) = ve_{21} + \alpha e_{12}$ for some $\alpha\in F$. Then $\varphi$ is an~automorphism on the algebra~$\mathbb{O}$.}

\medskip\textbf{Proposition 11.} \textit{Let $\varphi:\mathbb{O}\to \mathbb{O}$ be a~linear map such that $\varphi(e_{ii}) = e_{ii}$, $\varphi(ve_{ii}) = ve_{ii}$ for $i\in\{1,2\}$, $\varphi(e_{21}) = e_{21}$, $\varphi(ve_{21}) = ve_{21}$, $\varphi(e_{12}) = e_{12}+\alpha ve_{21}$, $\varphi(ve_{12}) = ve_{12} + \alpha e_{21}$ for some $\alpha\in F$. Then $\varphi$ is an~automorphism on the algebra~$\mathbb{O}$.}

\medskip\textbf{Proposition 12.} \textit{Let $\varphi:\mathbb{O}\to \mathbb{O}$ be a~linear map such that $\varphi(e_{ii}) = e_{ii}$ for $i\in\{1,2\}$, $\varphi(e_{12}) = e_{21}$, $\varphi(e_{21}) = e_{12}$, $\varphi(ve_{11}) = - ve_{21}$, $\varphi(ve_{12}) = - ve_{22}$, $\varphi(ve_{21}) = ve_{11}$, $\varphi(ve_{22}) = ve_{12}$ for some $\alpha\in F$. Then $\varphi$ is an~antiautomorphism on the algebra $\mathbb{O}$.}

\medskip\textbf{Proposition 13.} \textit{Let $\varphi:\mathbb{O}\to \mathbb{O}$ be a~linear map such that $\varphi(e_{ii}) = e_{ii}$, $\varphi(e_{12}) = e_{12} - \alpha ve_{22}$, $\varphi(e_{21}) = e_{21}$, $\varphi(e_{22}) = e_{22}$, $\varphi(ve_{11}) = ve_{11} + \alpha e_{21}$, $\varphi(ve_{12}) = ve_{12}$, $\varphi(ve_{21}) = ve_{21}$, $\varphi(ve_{22}) = ve_{22}$ for some $\alpha\in F$. Then $\varphi$ is an~automorphism on the algebra $\mathbb{O}$.}

\medskip\textbf{Proposition 14.} \textit{Let $\varphi:\mathbb{O}\to \mathbb{O}$ be a~linear map such that $\varphi(e_{ii}) = e_{ii}$, $\varphi(e_{12}) = e_{12}$, $\varphi(e_{21}) = e_{21} + \alpha ve_{11}$, $\varphi(e_{22}) = e_{22}$, $\varphi(ve_{11}) = ve_{11}$, $\varphi(ve_{12}) = ve_{12}$, $\varphi(ve_{21}) = ve_{21}$, $\varphi(ve_{22}) = ve_{22} - \alpha e_{12}$ for some $\alpha\in F$. Then $\varphi$ is an~automorphism on the algebra $\mathbb{O}$.}

\medskip\textbf{Proposition 15.} \textit{Let $\varphi:\mathbb{O}\to \mathbb{O}$ be a~linear map such that $\varphi(e_{11}) = e_{11} - \alpha ve_{12}$, $\varphi(e_{12}) = e_{12}$, $\varphi(e_{21}) = e_{21} + \alpha ve_{22}$, $\varphi(e_{22}) = e_{22} + \alpha ve_{12}$, $R(ve_{11}) = ve_{11} - \alpha e_{12}$, $R(ve_{12}) = ve_{12}$, $\varphi(ve_{21}) = ve_{21}-\alpha e_{11}+\alpha e_{22}+\alpha^2 ve_{12}$, $\varphi(ve_{22}) = ve_{22}$ for some $\alpha\in F$. Then $\varphi$ is an~automorphism on the algebra $\mathbb{O}$.
}

\medskip\textbf{Proposition 16.} \textit{Let $\varphi:\mathbb{O}\to \mathbb{O}$ be a~linear map such that $\varphi(e_{11}) = e_{11} + \alpha e_{12}$, $\varphi(e_{12}) = e_{12}$, $\varphi(e_{21}) = e_{21} - \alpha e_{11} + \alpha e_{22} - \alpha^2 e_{12}$, $\varphi(e_{22}) = e_{22} - \alpha e_{12}$, $R(ve_{11}) = ve_{11}$, $R(ve_{12}) = ve_{12}$, $\varphi(ve_{21}) = ve_{21}-\alpha ve_{11}$, $\varphi(ve_{22}) = ve_{22} - \alpha ve_{12}$ for some $\alpha\in F$. Then $\varphi$ is an~automorphism on the algebra $\mathbb{O}$.
}

\medskip\textbf{Proposition 17.} \textit{Let $\varphi:\mathbb{O}\to \mathbb{O}$ be a~linear map such that $\varphi(e_{11}) = e_{11} + \alpha e_{21}$, $\varphi(e_{12}) = e_{12} - \alpha e_{11} + \alpha e_{22} - \alpha^2 e_{21}$, $\varphi(e_{21}) = e_{21}$, $\varphi(e_{22}) = e_{22} - \alpha e_{21}$, $R(ve_{11}) = ve_{11} + \alpha ve_{21}$, $R(ve_{12}) = ve_{12} + \alpha ve_{22}$, $\varphi(ve_{21}) = ve_{21}$, $\varphi(ve_{22}) = ve_{22}$ for some $\alpha\in F$. Then $\varphi$ is an~automorphism on the algebra $\mathbb{O}$.}

\medskip If $A$ is an~algebra over $F$ and $R:A\to A$ is a~linear (over $F$) map, then $R$ is called \textit{\textbf{Rota-Baxter operator}} (\textit{\textbf{RB-operator}}) of a~weight $\lambda$ if for any $x,y\in A$ we have the following identity
\[R(x)R(y)=R(R(x)y+xR(y)+\lambda xy).\]

In this article we are interested only in RB-operators of zero weight, $\lambda = 0$.

An image of $R$ is a~subalgebra $B$ of an algebra $A$ and $\ker (R)$ is a~$\mathrm{Im}(R)$-bimodule.

\medskip In \cite{BGP} RB-operators on the algebra $M_2(F)$ for algebraically closed $F$ were described. We will need this description for any field with characteristic $\neq 2$.

\medskip\textbf{Proposition 18.} \textit{Let $R: M_2(F) \to M_2(F)$ be a~Rota-Baxter operator on $M_2(F)$ of weight zero. Then, up to conjugation by automorphism, antiautomorphism and up to multiplication by a~scalar, $R$ acts in one of the following ways:
\begin{enumerate}
    \item $R(e_{21})=e_{11}, \quad R(e_{11})=R(e_{12})=R(e_{22})=0$;
    \item $R(e_{21})=e_{12}, \quad R(e_{11})=R(e_{12})=R(e_{22})=0$;
    \item $R(e_{21})=e_{11}, \quad R(e_{22})=e_{12}, \quad R(e_{11})=R(e_{12})=0$;
    \item $R(e_{21})=-e_{11}, \quad R(e_{11})=e_{12}, \quad R(e_{12})=R(e_{22})=0$.
\end{enumerate}
}

\textbf{Proof.} In (\cite{G2017}, Lemma 1), it is shown that on a~simple unital finite-dimensional non-one-dimensional algebra the kernel of any Rota-Baxter operator of weight 0 has dimension at least two, moreover, the unit does not lie in the image of this operator. Thus, by the theorem on the dimension of the kernel and the image of the linear map, $\dim(\mathrm{Im}(R))\le 2$, and $\mathrm{Im}(R)$ does not contain the identity matrix.

The results of paper \cite{Octonions} imply that any one-dimensional non-unital subalgebra in $\mathbb{O}$ has the form $Fe_{11}$ or $Fe_{12}$ up to automorphism and antiautomorphism. Thus, subalgebras in $M_2(F)$ also have the form $Fe_{11}$ or $Fe_{12}$ up to automorphism and antiautomorphism. The results of paper \cite{Octonions} imply that, up to automorphism and antiautomorphism, any two-dimensional non-unital subalgebra in $\mathbb{O}$ either has the form $Fe_{11}+Fe_{12}$ or has zero multiplication. However, $M_2(F)$ does not have two-dimensional subalgebras with zero multiplication. Thus, the only two-dimensional subalgebra in $M_2(F)$ (up to automorphism and antiautomorphism) has the form $Fe_{11}+Fe_{12}$.

Since $\mathrm{Im}(R)$ is a~subalgebra, the following cases are possible.

\smallskip 1) $\mathrm{Im}(R)=Fe_{11}$. Let $R(e_{ij})=\alpha_{ij} e_{11}$ for any $i,j$. Then
\[\alpha_{11}^2 e_{11} = R(e_{11})R(e_{11}) = R(R(e_{11})e_{11}+e_{11}R(e_{11}))= 2\alpha^2 e_{11},\]
whence $\alpha_{11}=0$. Further,
\[\alpha_{22}^2 e_{11} = R(e_{22})R(e_{22}) = R(R(e_{22})e_{22}+e_{22}R(e_{22})) = 0,\]
whence $\alpha_{22}=0$. Further
\[\alpha_{12}\alpha_{21} e_{11} = R(e_{12})R(e_{21}) = R(R(e_{12})e_{21}+e_{12}R(e_{21})) = 0,\]
whence $\alpha_{12}\alpha_{21}=0$. Up to antiautomorphism (transposition), we can assume that $\alpha_{12}=0$. Up to multiplication by a~scalar, we can assume that $R(e_{21})=e_{11}$.

\smallskip 2) $\mathrm{Im}(R)=Fe_{12}$. Let $R(e_{ij})=\alpha_{ij} e_{12}$ for any $i,j$. Then
\[0 = R(e_{12})R(e_{22}) = R(\alpha_{12}e_{12}) = \alpha_{12}^2 e_{12},\]
whence $\alpha_{12}=0$. Further,
\[0 = R(e_{11})R(e_{21}) = R(\alpha_{11}e_{11}+\alpha_{21}e_{12}) = \alpha_{11}^2 e_{12},\]
whence $\alpha_{11}=0$. Further,
\[0 = R(e_{21})R(e_{22}) = R(\alpha_{21}e_{12} + \alpha_{22}e_{22}) = \alpha_{22}^2 e_{12},\]
whence $\alpha_{22}=0$. 

Up to multiplication by a~scalar, we can assume that $R(e_{21}) = e_{12}$.

\smallskip 3) $\mathrm{Im}(R) = Fe_{11}+Fe_{12}$. Let $x=\alpha e_{11}+\beta e_{12}+\gamma e_{21}+\delta e_{22}\in\ker (R)$. Since $e_{11},e_{12}\in\mathrm{Im}(R)$ and $\ker (R)$ is a~$\mathrm{Im}(R)$-bimodule, then $e_{11}x=\alpha e_{11}+\beta e_{12}\in\ker (R)$, whence $e_{11}xe_{11}=\alpha e_{11}\in\ker (R)$, so $e_{11}x-e_{11}xe_{11}=\beta e_{12}\in\ker (R)$. Similarly $\alpha e_{11} + \gamma e_{21}\in\ker (R)$ and $\gamma e_{21}\in\ker (R)$, whence $\delta e_{22}\in\mathrm{Ker(R)}$. Thus, if $\alpha e_{11}+\beta e_{12}+\gamma e_{21}+\delta e_{22}\in\ker (R)$, then $\alpha e_{11},\beta e_{12},\gamma e_{12}, \delta e_{22}\in\ker (R)$.

Suppose that $e_{11}\in\ker (R)$. Then $e_{12}=e_{11}\cdot e_{12}\in\ker (R)$ and $\ker (R)=Fe_{11}+Fe_{12}$. Let $R(e_{22})=\xi_4 e_{11}+\eta_4 e_{12}$, $R(e_{21})=\xi_3 e_{11}+\eta_3 e_{12}$. Then 
\[\xi_4^2 e_{11}+\xi_4\eta_4 e_{12} = R(e_{22})R(e_{22}) = R(\eta_4 e_{12}) = 0,\]
whence $\xi_4=0$. Further, 
\[\xi_3^2 e_{11}+\xi_3\eta_3 e_{12} = R(e_{21})R(e_{21}) = R(\eta_3 e_{11} + \xi_3 e_{21} + \eta_3 e_{22}) = \xi_3^2 e_{11} + (\xi_3\eta_3+\eta_3\eta_4)e_{12},\] 
whence $\eta_3\eta_4=0$. Finally,
\[\xi_3\eta_4 e_{12} = R(e_{21})R(e_{22}) = R(\eta_3e_{12} + \eta_4e_{22}) = \eta_4^2 e_{12},\]
whence $\eta_4(\eta_4-\xi_3)=0$. Since $\dim (\mathrm{Im}(R)) = 2$, then $\eta_4\neq 0$ and, up to multiplication by a~scalar, $R(e_{21})=e_{11}$, $R(e_{22})=e_{12}$.

Suppose that $e_{11}\notin\ker (R)$. Then $e_{21}\notin\ker (R)$ (otherwise $e_{11}=e_{12}\cdot e_{21}\in\ker (R)$), so $\ker (R)=Fe_{12}+Fe_{22}$. Let $R(e_{11})=\xi_1 e_{11} + \eta_1 e_{12}$, $R(e_{21})=\xi_3 e_{11} + \eta_3 e_{12}$. Then
\[\xi_1^2 e_{11} + \xi_1\eta_1 e_{12} = R(e_{11})R(e_{11}) = R(2\xi_1 e_{11} + \eta_1 e_{12}) = 2\xi_1^2 e_{11} + 2\xi_1\eta_1 e_{12},\]
whence $\xi_1=0$. Further, 
\[\xi_3^2 e_{11}+\xi_3\eta_3 e_{12} = R(e_{21})R(e_{21}) = R(\eta_3 e_{11} + \xi_3 e_{21} + \eta_3 e_{22}) = \xi_3^2 e_{11} + (\eta_3\eta_1 + \xi_3\eta_3)e_{12},\] 
whence $\eta_3\eta_1 = 0$. Since $\mathrm{dim} V = 2$, then $\eta_1\neq 0$ and $\eta_3=0$. Finally,
\[0=R(e_{11})R(e_{21}) = R(\eta_1 e_{11} + \xi_3 e_{11}) = (\eta_1+\xi_3)\eta_1 e_{12}.\]
Since $\eta_1\neq 0$, then, up to multiplication by a~scalar, $R(e_{11})=e_{12}$, $R(e_{21})=-e_{11}$. The proposition is proven.

\section{RB-Operators with one-dimensional image}

In \cite{Octonions} it was proved that there are only two one-dimensional non-unital subalgebras $B$ in $\mathbb{O}$, up to action of automorphism: nilpotent $Fe_{12}$ and idempotent $Fe_{11}$. Let us describe the Rota-Baxter operators of zero weight on $\mathbb{O}$ with these images.

\medskip\textbf{Lemma 1.} \textit{Let $R$ be a~Rota-Baxter operator of zero weight on the split Cayley-Dickson algebra $\mathbb{O}$ and $\mathrm{Im}(R)=Fe_{12}$. Then, up to conjugation by automorphism, antiautomorphism and up to multiplication by a~scalar, an operator $R$ acts in one of the following ways (an operator $R$ is zero on unspecified basic elements $e_{ij}$, $ve_{ij}$):
\begin{enumerate}
    \item $R(e_{21}) = e_{12}$;
    \item $R(ve_{22}) = e_{12}$.
\end{enumerate}
}

\textbf{Proof.} Since $\mathrm{Im}(R)\subset M_2(F)$, then $R|_{M_2(F)}$ is the Rota-Baxter operator on the subalgebra $M_2(F)$. According to the Proposition~18, we can assume that $R(e_{11})=R(e_{12})=R(e_{22})=0$.

Note that $(\mathrm{Im}(R))^2=0$, so $e_{12}R(y)=R(x)e_{12}=R(x)R(y)= 0$ for any $x,y\in\mathbb{O}$. Let $R(ve_{12})=\alpha e_{12}$, then
\[0=R(ve_{22})R(ve_{12})=R(ve_{22}\cdot \alpha e_{12})=\alpha R(ve_{12})=\alpha^2 e_{12},\]
whence $\alpha = 0$ and $R(ve_{12})=0$. Let $R(ve_{11})=\beta e_{12}$, then
\[0=R(ve_{11})R(ve_{21})=R(\beta e_{12} \cdot ve_{21} + ve_{11}R(ve_{21})) = -\beta R(ve_{11})=-\beta^2 e_{12},\]
whence $\beta=0$ and $R(ve_{11})=0$. Let $R(ve_{21})=\alpha_1 e_{12}$ and $R(ve_{22})=\alpha_2 e_{12}$. 

If $\alpha_2\neq 0$, then there exists $\varepsilon\in F$ such that $ve_{21}+\varepsilon ve_{22}\in \ker (R)$. According to the Proposition~2, we can assume that $\varepsilon=0$, that is, $\alpha_1 = 0$. Thus, we can assume that either $\alpha_2=0$ or $\alpha_1=0$. According to the Proposition~4, we can assume that $\alpha_1=0$. Thus, $R(ve_{22})=\alpha_2 e_{12}$ and $R(e_{21}) = \alpha_3 e_{12}$.

Let $\alpha_2 = 0$. Then we can assume that $R(e_{21}) = e_{12}$.

Let $\alpha_3 = 0$. Then we can assume that $R(ve_{22}) = e_{12}$. 

Let $\alpha_2\neq 0$ and $\alpha_3\neq 0$. Up to multiplication by $\frac{1}{\alpha_2}$, we can assume that $R(ve_{22}) = e_{12}$, $R(e_{21}) = \alpha_3 e_{12}$. Conjugation by automorphism from Proposition~6 with $\alpha = \frac{1}{\alpha_3}$ gives us $R(ve_{22}) = \frac{1}{\alpha_3} e_{12}$, $R( e_{21}) = \frac{1}{\alpha_3} e_{12}$, which after multiplication by $\alpha_3$ gives $R(ve_{22})=e_{12}$, $R(e_{ 21})=e_{12}$. Proposition~15 with a~scalar $\alpha = 1$ allows us to assume that $R(e_{21}) = 0$ and $R(ve_{22}) = e_{12}$. The lemma is proven.

\medskip\textbf{Lemma 2.} \textit{Let $R$ be a~Rota-Baxter operator of zero weight on the split Cayley-Dickson algebra $\mathbb{O}$ and $\mathrm{Im}(R)=Fe_{11}$. Then, up to conjugation by automorphism, antiautomorphism and up to multiplication by a~scalar, an operator $R$ acts in the following way (an operator $R$ is zero on unspecified basic elements $e_{ij}$, $ve_{ij}$):
\[ R(e_{21}) = e_{11}.\]
} 

\textbf{Proof.} Since $\mathrm{Im}(R)\subset M_2(F)$, then $R$ is the Rota-Baxter operator on the subalgebra $M_2(F)$. According to the Proposition~18, we can assume that $R(e_{11})=R(e_{12})=R(e_{22})=0$. Let us introduce the notation $R(ve_{ij})=\beta_{ij} e_{11}$. Then
\[\beta_{12}\beta_{22}e_{11}=R(ve_{12})R(ve_{22})=R(\beta_{12} ve_{22} + \beta_{22} ve_{12}) = 2\beta_{12}\beta_{22}e_{11},\]
whence $\beta_{12}\beta_{22}=0$. Further,
\[\beta_{12}\beta_{21}e_{11}=R(ve_{21})R(ve_{12})=R(0)=0,\]
whence $\beta_{12}\beta_{21}=0$. Further,
\[\beta_{11}\beta_{22}e_{11}=R(ve_{22})R(ve_{11})=R(0)=0,\]
whence $\beta_{11}\beta_{22}=0$. Further,
\[\beta_{11}\beta_{21}e_{11}=R(ve_{11})R(ve_{21})=R(\beta_{11}ve_{21} + \beta_{21}ve_{11})=2\beta_{11}\beta_{21}e_{11},\]
whence $\beta_{11}\beta_{21}=0$. Let $R(e_{21})=\alpha e_{11}$. Then for any $j\in\{1,2\}$ we have
\[\alpha\beta_{2j}e_{11} = R(ve_{2j})R(e_{21}) = 0,\]
whence $\alpha\beta_{21}=\alpha\beta_{22} = 0$. Thus, either $R(e_{21})=R(ve_{11})=R(ve_{12})=0$ or $R(ve_{21})=R(ve_{22})=0$. The involution from Proposition~1 allows us to assume that $R(ve_{21})=R(ve_{22})=0$.

If $\beta_{12}\neq 0$, then there exists $\varepsilon_1\in F$ such that $ve_{11}+\varepsilon_1 ve_{12}\in\ker (R)$. By Proposition~2 we can assume that $\varepsilon_1 = 0$, so $ve_{11}\in\ker (R)$. Thus, in any case, either $\beta_{11}=0$ or $\beta_{12}=0$. By Proposition~4 we can assume that $\beta_{12}=0$. Conjugation by automorphism from Proposition~6, as above, allows us to assume that either $R(e_{21})=e_{11}$ and $R(ve_{11}) = e_{11}$, or $R(ve_{11})=0$ and $R(e_{21})=e_{11}$, or $R(e_{21})=0$ and $R(ve_{11})=e_{11}$. The composition of the classical involution and the involution from Proposition~5 allows us to consider that the second and third cases are equivalent. The Proposition~13 with a~scalar $1$ states that the first and the second cases are equivalent. The lemma is proven.

\medskip

\section{RB-Operators with two-dimensional image}

In \cite{Octonions} it was proved that there are only two two-dimensional non-unital subalgebras $B$ in $\mathbb{O}$, up to action of automorphism and antiautomorphism: idempotent $Fe_{11}+Fe_{12}$ and nilpotent $Fve_{12}+Fve_{22}$. Let us describe the Rota-Baxter operators of zero weight on $\mathbb{O}$ with these images.

\medskip\textbf{Lemma 3.} \textit{Let $R$ be a~Rota-Baxter operator of zero weight on the split Cayley-Dickson algebra $\mathbb{O}$ and $\mathrm{Im}(R)=Fe_{11}+Fe_{12}$. Then, up to conjugation by automorphism, antiautomorphism and up to multiplication by a~scalar, an operator $R$ acts in one of the following ways for some $\alpha\in F$ (an operator $R$ is zero on unspecified basic elements $e_{ij}$, $ve_{ij}$):
\begin{enumerate}
    \item $R(e_{21})=e_{11}$, $R(e_{22})=e_{12}$;
    \item $R(e_{21})=-e_{11}$, $R(e_{11})=e_{12}$;
    \item $R(e_{21})=e_{11}$, $R(ve_{21})= e_{12}$;
    \item $R(ve_{11})=\alpha e_{11}$, $R(ve_{21})= e_{12}$, $\alpha\neq 0$;
    \item $R(ve_{11}) = e_{12}, R(ve_{21}) = e_{11}$;
    \item $R(ve_{21}) = \alpha e_{11}, R(ve_{22}) = e_{12}$, $\alpha\neq 0$.
\end{enumerate}
}

\textbf{Proof.} According to the Proposition~18, the following cases are possible.

1) $R(e_{21})=e_{11}$, $R(e_{22})=e_{12}$, $R(e_{11})=R(e_{12})=0$. 

Let $R(ve_{11})=\alpha_1 e_{11} + \beta_1 e_{12}$. Then
\[\alpha_1 e_{11} = R(ve_{11})R(e_{21}) = R(\beta_1 e_{11} + ve_{11}) = \alpha_1 e_{11} + \beta_1 e_{12},\]
whence $\beta_1 = 0$. Further,
\[\alpha_1 e_{12} = R(ve_{11})R(e_{22}) = R(0) = 0,\]
whence $\alpha_1 = 0$. Therefore, $R(ve_{11})=0$.

Let $R(ve_{12})=\alpha_2 e_{11} + \beta_2 e_{12}$. Then
\[\alpha_2 e_{11} = R(ve_{12})R(e_{21}) = R(\beta_2 e_{11} + ve_{12}) = \alpha_2 e_{11} + \beta_2 e_{12},\]
hence $\beta_2 = 0$. Further,
\[\alpha_2 e_{12} = R(ve_{12})R(e_{22}) = R(0) = 0,\]
hence $\alpha_2 = 0$. Therefore, $R(ve_{12})=0$.

Let $R(ve_{21})=\alpha_3 e_{11} + \beta_3 e_{12}$. Then
\[\alpha_3 e_{11} = R(ve_{21})R(e_{21}) = R(\beta_3 e_{11}) = 0,\]
hence $\alpha_3 = 0$. Further,
\[\beta_3 e_{12} = R(e_{21})R(ve_{21}) = R(\beta_3 e_{22} + ve_{21}) = 2\beta_3 e_{12},\]
hence $\beta_3 = 0$. Therefore, $R(ve_{21})=0$.

Let $R(ve_{22})=\alpha_4 e_{11} + \beta_4 e_{12}$. Then
\[\alpha_4 e_{11} = R(ve_{22})R(e_{21}) = R(\beta_4 e_{11}) = 0,\]
hence $\alpha_4=0$. Further,
\[\beta_4 e_{12} = R(e_{21})R(ve_{22}) = R(ve_{22}+\beta_4 e_{22}) = 2\beta_4 e_{12},\]
hence $\beta_4=0$. Therefore, $R(ve_{22})=0$. We have obtained case (1) from the statement of the lemma.

\smallskip 2) $R(e_{21})=-e_{11}$, $R(e_{11})=e_{12}$, $R(e_{12})=R(e_{22})=0$.

Let $R(ve_{11})=\alpha_1 e_{11} + \beta_1e_{12}$. Then
\[0 = R(e_{11})R(ve_{11}) = R(e_{11} R(ve_{11})) = R(\alpha_1 e_{11} + \beta_1 e_{12}) = \alpha_1 e_{12},\]
hence $\alpha_1 = 0$. Further, 
\[ - \beta_1 e_{12} = R(e_{21}) R(ve_{11}) = R(\beta_1 e_{22}) = 0,\]
hence $\beta_1 = 0$. Therefore, $R(ve_{11})=0$. 

Let $R(ve_{22}) = \alpha_4 e_{11} + \beta_4 e_{22}$. Then
\[-\alpha_4 e_{11} - \beta_4 e_{12} = R(e_{21})R(ve_{22}) = -2\alpha_4 e_{11} - \beta_4 e_{12},\]
hence $\alpha_4 = 0$. Further,
\[0 = R(ve_{22})R(e_{21}) = R(\beta_4 e_{11} - ve_{22}\cdot e_{11}) = \beta_4 e_{12},\]
hence $\beta_4 = 0$. Therefore, $R(ve_{22})=0$.

Let $R(ve_{12}) = \alpha_2 e_{11} + \beta_2 e_{12}$. Then
\[ 0 = R(e_{11})R(ve_{12}) = R(e_{12}\cdot ve_{12} + \alpha_2 e_{11} + \beta_2 e_{12}) = \alpha_2 e_{12},\]
hence $\alpha_2 = 0$. Further,
\[0 = R(e_{11})R(ve_{22}) = R(-ve_{12} + \beta_4 e_{12}) = -\beta_2 e_{12},\]
hence $\beta_2 = 0$. Therefore, $R(ve_{12})=0$. 

Let $R(ve_{21}) = \alpha_3 e_{11} + \beta_3 e_{12}$. Then
\[0 = R(e_{11})R(ve_{21}) = R(e_{12}\cdot ve_{21} + \alpha_3 e_{11} + \beta_3 e_{12}) = \alpha_3 e_{12},\]
hence $\alpha_3 = 0$. Further,
\[0 = R(ve_{21})R(e_{21}) = R(\beta_3 e_{11} - ve_{21}\cdot e_{11}) = \beta_3 e_{12},\]
hence $\beta_3 = 0$. Therefore, $R(ve_{21})=0$. We have obtained case (2) from the statement of the lemma.

\smallskip In the remaining three cases, the dimension of $R(M_2(F))$ does not exceed one. Then the dimension of $R(vM_2(F))$ is not less than one. Let $0\neq x = \alpha ve_{11}+\beta ve_{12}+\gamma ve_{21}+\delta ve_{22}\in\ker (R)$. Since $e_{11},e_{12}\in\mathrm{Im}(R)$ and $\ker (R)$ is an $\mathrm{Im}(R)$-bimodule, then $e_ {11}x=\gamma ve_{21} + \delta ve_{22}\in\ker (R)$, whence $\alpha ve_{11} + \beta ve_{12}\in\mathrm{ Ker}(R)$. But then $e_{12}(e_{11}x) = -\gamma ve_{11}-\delta ve_{12}\in\ker (R)$. Let $V_1=\ker (R)\cap (Fve_{11}+Fve_{12})$, $V_2=\ker (R)\cap (Fve_{21}+Fve_{22})$. Thus, $\ker (R)\cap vM_2(F) = V_1\oplus V_2$, and $\dim V_1 \ge \dim V_2$. Since $2\le \dim (\ker (R)\cap vM_2(F)) \le 3$ (by the theorem on the dimension of the kernel and image for $R|_{vM_2(F)}$), then either $\dim V_1 = 2$ and $\dim V_2 = 1$, or $\dim V_1=2$ and $\dim V_2 = 0$, or $\dim V_1 = \dim V_2 = 1$. In the first case, $ve_{11},ve_{12}\in\ker (R)$ and we can assume (by Propositions~2--4) that $ve_{21}\in\mathrm{Ker }(R)$. In the second case, $ve_{11},ve_{12}\in\ker (R)$. In the third case, we can assume (by Propositions~2--4) that $ve_{11},ve_{21}\in\ker (R)$.

\smallskip 3) $R(e_{21})=e_{11}$, $R(e_{11})=R(e_{12})=R(e_{22})=0$. 

Let $R(ve_{11})=\alpha_1 e_{11} + \beta_1 e_{12}$. Then
\[\alpha_1 e_{11} = R(ve_{11})R(e_{21}) = R(\beta_1 e_{11} + ve_{11})=R(ve_{11}) = \alpha_1 e_{11} + \beta_1 e_{12},\]
hence $\beta_1 = 0$. Let $R(ve_{12})=\alpha_2 e_{11} + \beta_2 e_{12}$. Then
\[\alpha_2 e_{11} = R(ve_{12})R(e_{21}) = R(ve_{12}R(e_{21})) = R(ve_{12}) = \alpha_2 e_{11} + \beta_2 e_{12},\]
hence $\beta_2 = 0$. Let $R(ve_{21})=\alpha_3 e_{11} + \beta_3 e_{12}$. Then
\[ \alpha_3 e_{11} = R(ve_{21}) R(e_{21}) = R(\beta_3 e_{11}) = 0,\]
hence $\alpha_3 = 0$. Let $R(ve_{22})=\alpha_4 e_{11} + \beta_4 e_{12}$. Then
\[\alpha_4 e_{11} = R(ve_{22})R(e_{21}) = R(\beta_4 e_{11}) = 0,\]
hence $\alpha_4=0$. Further,
\[0=R(ve_{22})R(ve_{21})=R(-\beta_4 ve_{11} + \beta_3 ve_{12}) = (\beta_3\alpha_2 -\beta_4 \alpha_1)e_{11},\]
hence $\beta_3\alpha_2=\beta_4\alpha_1$. 

3a) Let $\dim V_1 = 2$ and $\dim V_2 = 1$. Then, by above and by Proposition~4, we can assume that
\begin{gather*}
R(e_{11})=R(e_{12})=R(e_{22})=R(ve_{11})=R(ve_{12})=R(ve_{22})=0,\\
R(e_{21})=e_{11}, \quad R(ve_{21})=\beta_4 e_{12}.
\end{gather*}
Consider the automorphism $\varphi$ from Proposition~6 for $\alpha=\frac{1}{\beta_4}$. Then $\varphi^{-1}R\varphi(e_{21})=\frac{1}\beta_4 e_{11}$, $\varphi^{-1}R\varphi(ve_{21}) = \frac{1}{\beta_4} e_{12}$. After multiplication by the scalar $\beta_4$ we can assume that
\[R(e_{21})=e_{11}, \quad R(ve_{21})= e_{12}.\]
We have obtained case (3) in the statement of the lemma.

3b) Let $\dim V_1 = 2$ and $\dim V_2 = 0$. If we consider the restriction of the mapping $R$ to $V_2$, we will find that the image of this mapping has a~dimension at most one, which implies (according to the the theorem on the dimension of the kernel and the image) that the kernel must have a~dimension at least one. Therefore, this case is impossible.

3c) Let $\dim V_1 = \dim V_2 = 1$.  Then, by above, we can assume that
\begin{gather*}R(e_{11})=R(e_{12})=R(e_{22})=R(ve_{11})=R(ve_{22})=0,\\
R(e_{21})=e_{11},\quad R(ve_{11})=\alpha_2 e_{11},\quad R(ve_{21})=\beta_4 e_{12}.
\end{gather*}
Consider the automorphism $\varphi$ from Proposition~6 for $\alpha=\frac{1}{\beta_4}$. Then $\varphi^{-1}R\varphi(e_{21})=\frac{1}\beta_4 e_{11}$, $\varphi^{-1}R\varphi(ve_{11}) = \alpha_2 e_{11}$, $\varphi^{-1}R\varphi(ve_{21}) = \frac{1}{\beta_4} e_{12}$. After multiplication by the scalar $\beta_4$ we can assume that (here $\varepsilon = \alpha_2\beta_4$)
\[R(e_{21})=e_{11}, \quad R(ve_{11}) = \varepsilon e_{11}, \quad R(ve_{21})=e_{12}.\]
The Proposition~14 with a~scalar $\varepsilon$ allows us to assume that $R(e_{21}) = 0$. 
We have obtained the case (4) in the statement of the lemma.

\smallskip 4) $R(e_{21})=e_{12}$, $R(e_{11}) = R(e_{22}) = R(e_{12}) = 0$. Let $R(ve_{11})=\alpha_1 e_{11} + \beta_1 e_{12}$. Then
\[\alpha_1 e_{12} = R(ve_{11})R(e_{21}) = R(\beta_1 e_{11}) = 0,\]
whence $\alpha_1=0$. Let $R(ve_{12})=\alpha_2 e_{11} + \beta_2 e_{12}$. Then
\[\alpha_2 e_{12} = R(ve_{12})R(e_{21}) = R(\beta_2 e_{11}) = 0,\]
whence $\alpha_2=0$. Let $R(ve_{21})=\alpha_3 e_{11} + \beta_3 e_{12}$. Then
\[\alpha_3 e_{12} = R(ve_{21})R(e_{21}) = R(\beta_3 e_{11} + ve_{11}) = \beta_1 e_{121},\]
whence $\alpha_3=\beta_1$. Let $R(ve_{22})=\alpha_4 e_{11} + \beta_4 e_{12}$. Then
\[\alpha_4 e_{12} = R(ve_{22})R(e_{21}) = R(\beta_4 e_{11} + ve_{12}) = \beta_2 e_{12},\]
whence $\alpha_4 = \beta_2$.

4a) Let $\dim V_1 = 2$. By above we have $\beta_1=\beta_2=0$, hence $\alpha_3=\alpha_4=0$. Then $\mathrm{Im}(R)=Fe_{12}$, a~contradiction. Therefore, this case is impossible.

4b) Let $\dim V_1 = \dim V_2 = 1$. By above and Proposition~4, we can assume that $ve_{12}, ve_{22}\in\ker (R)$. Therefore, $\beta_2=\alpha_4=\beta_4=0$. We have
\begin{gather*}
R(e_{11})=R(e_{12})=R(e_{22})=R(ve_{12})=R(ve_{22})=0,\\
R(e_{21})=e_{12},\quad R(ve_{11})=\beta_1 e_{12},\quad R(ve_{21})=\beta_1 e_{11} + \beta_3 e_{12}.
\end{gather*}
Consider the automorphism $\varphi$ from Proposition~6 for $\alpha=\beta_1$. Then $\varphi^{-1}R\varphi(e_{21}) = \beta_1^2 e_{12}$, $\varphi^{-1}R\varphi(ve_{11}) = \beta_1^ 2 e_{12}$, $R(ve_{21}) = \beta_1^2 e_{11} + \beta_3\beta_1^2 e_{12}$. After multiplication by $\frac{1}{\beta_1^2}$ we can assume that
\[R(e_{21}) = e_{12}, \quad R(ve_{11}) = e_{12}, \quad R(ve_{21}) = e_{11} + \beta_3 e_{12}.\]
Proposition~14 allows us to assume that
\[R(e_{21}) = 0, \quad R(ve_{11}) = e_{12}, \quad R(ve_{21}) = e_{11} + \beta_3 e_{12}.\]
Proposition~16 with a~scalar $\alpha = -\frac{\beta_3}{2}$ allows us to assume that
\[R(e_{21}) = 0, \quad R(ve_{11}) = e_{12}, \quad R(ve_{21}) = e_{11}.\]
We have obtained case (5) in the statement of the lemma.

\smallskip 5) $R(M_2(F))=0$. By above we can assume that either $\ker (R)\cap (vM_2(F)) = Fve_{11} + Fve_{12}$, or $\ker (R)\cap (vM_2(F)) = Fve_{11} + Fve_{21}$. Let $R(ve_{12})=\alpha_2 e_{11} + \beta_2 e_{12}$, $R(ve_{21})=\alpha_3 e_{11} + \beta_3 e_{12}$, $R(ve_{22})=\alpha_4 e_{11} + \beta_4 e_{12}$. 

5a) $R(ve_{21})=0$, that is $\alpha_3=\beta_3=0$. Then
\[
\alpha_2\alpha_4 e_{11} + \alpha_2\beta_4 e_{12} =R(ve_{12})R(ve_{22}) = (\alpha_2\alpha_4 -\beta_2\alpha_2 + \alpha_4\alpha_2)e_{11} + (\alpha_2\beta_4 - \beta_2\beta_2 + \alpha_4 \beta_2)e_{12}.
\]
Since $R(ve_{12})=\alpha_2 e_{11} + \beta_2 e_{12}\neq 0$, then $\alpha_4=\beta_2$.

If $\beta_2 = 0$, then we can assume (up to Proposition~4) that $R(ve_{11})=\alpha_2 e_{11}$, $R(ve_{21}) = e_{12 }$. We have obtained case (4) from the statement of the lemma.

If $\beta_2\neq 0$, then conjugation by automorphism from Proposition~7 with $\alpha = \frac{1}{\beta_2}$ gives us (with $\varepsilon = \frac{\gamma}{\beta^2}$)
\[R(ve_{12}) = \alpha_2 e_{11} + e_{12}, \quad R(ve_{22}) = e_{11} + \varepsilon e_{12}.\]
The conjugation by automorphism from Proposition~6 for $\alpha = \frac{1}{\varepsilon}$ and Proposition~4 allow us to assume that $R$ has the form
\[R(ve_{11}) = \alpha_2\varepsilon e_{11} + e_{12}, \quad R(ve_{21}) = e_{11} + e_{12}.\]
After conjugation with an automorphism from Proposition~16 with a~scalar $-\frac{1}{\alpha_2}$ and a~multiplication by a~scalar, we obtain case (4) in the statement of the lemma (if $\alpha_2\neq 1$, otherwise we obtain $R(ve_{21})=0$, it is a~contradiction).

It is easy to see that the operator $R$ with these conditions is a~Rota-Baxter operator. It remains to note that in order for the condition $\mathrm{Im}(R)=Fe_{11}+Fe_{12}$ to be satisfied, it is necessary and sufficient that $\alpha_2\varepsilon\neq 1$.

5b) $R(ve_{12})=0$, that is $\alpha_2=\beta_2=0$. Then
\begin{gather*}R(ve_{21}) = \alpha_3 e_{11} + \beta_3 e_{12},\\
R(ve_{22}) = \alpha_4 e_{11} + \beta_4 e_{12}.
\end{gather*}
It is easy to see that the operator $R$ with these conditions is a~Rota-Baxter operator. Note that in order for the condition $\mathrm{Im}(R)=Fe_{11}+Fe_{12}$ to be satisfied, it is necessary and sufficient that $\alpha_3\beta_4\neq \alpha_4\beta_3$. Next, let $\beta_3 = 0$. Then we can assume that
\begin{equation}\label{beta3=0}
    R(ve_{21})=e_{11}, \quad R(ve_{22}) = \alpha_4e_{11} + \beta_4 e_{12}.
\end{equation}
By the Proposition~7 we can assume that
\[R(ve_{21})=e_{11}, \quad R(ve_{22}) = e_{11} + \frac{\beta_4}{\alpha_4}e_{12}.\]
Then, by the Proposition~8, we can assume that
\[R(ve_{21})=\gamma e_{11}, R(ve_{22}) = e_{11} + e_{12},\]
where $\gamma = \frac{\beta_4}{\alpha_4}$. According to the Proposition~3 we can assume that
\[R(ve_{21}) = \gamma e_{11}, \quad R(ve_{22}+\gamma^{-1}ve_{21}) = e_{11} + e_{12},\]
where
\[R(ve_{21}) = \gamma e_{11}, \quad R(ve_{22}) = e_{12}.\]
We obtain case (6) in the statement of the lemma. Let $\beta_3 \neq 0$. Then we can assume that
\[R(ve_{21}) = \alpha_3 e_{11} + e_{12}.\]
Further, by the Proposition~3 we can assume that
\[R(ve_{21}) = \alpha_3 e_{11} + \beta_4 e_{12}, R(ve_{22} + \beta_4 ve_{21}) = \alpha_4 e_{11} + \beta_4 e_{12},\]
whence
\[R(ve_{21}) = \alpha_3 e_{11} + e_{12}, \quad R(ve_{22}) = \gamma e_{12},\]
where $\gamma = \alpha_4-\beta_4\alpha_3$. Then we can assume that
\[R(ve_{21}) = \alpha_3' e_{11} + \beta_3' e_{12}, \quad R(ve_{22}) = e_{11}.\]
According to the Proposition~4 we can assume that
\[R(ve_{21}) = e_{11}, \quad R(ve_{22}) = -\alpha_3' e_{11} - \beta_3' e_{12}.\]
This is exactly the already discussed case $\beta_3 = 0$, formula \eqref{beta3=0}. The lemma is proven.

\medskip\textbf{Corollary 1.} \textit{Let $R$ be a~Rota-Baxter operator of zero weight on the split Cayley-Dickson algebra $\mathbb{O}$ and $\mathrm{Im}(R)=Fe_{11}+Fe_{12}$. If a~field $F$ is quadratically closed then, up to conjugation by automorphism, antiautomorphism and up to multiplication by a~scalar, an operator $R$ acts in one of the following ways for some $\alpha\in F$ (an operator $R$ is zero on unspecified basic elements $e_{ij}$, $ve_{ij}$):
\begin{enumerate}
    \item $R(e_{21})=e_{11}$, $R(e_{22})=e_{12}$;
    \item $R(e_{21})=-e_{11}$, $R(e_{11})=e_{12}$;
    \item $R(e_{21})=e_{11}$, $R(ve_{21})= e_{12}$;
    \item $R(ve_{11})=e_{11}$, $R(ve_{21})= e_{12}$;
    \item $R(ve_{11})=e_{12}, R(ve_{21})= e_{11}$;
    \item $R(ve_{21}) =  e_{11}, R(ve_{22}) = e_{12}$;
\end{enumerate}
}
\textbf{Proof.} Let us consider the resulting cases on $R$ in Lemma~3. Cases 1--3,5--7 remained the same.

4) Proposition~7 with a~scalar $\sqrt{\alpha}$ allows us to assume that
\[R(ve_{11}) = \sqrt{\alpha}e_{11}, R(ve_{21}) = \sqrt {\alpha}e_{12}.\]
Up to multiplication by a~scalar, we obtain the required operator.

6) Proposition~7 with a~scalar $\sqrt{\alpha}$ allows us to assume that
\[R(ve_{21}) = \alpha e_{11}, R(ve_{22}) = \alpha e_{12}.\]
Up to multiplication by a~scalar, we obtain the required operator. The corollary is proven.

\medskip \textbf{Lemma 4.} \textit{Let $R$ be a~Rota-Baxter operator of zero weight on the split Cayley-Dickson algebra $\mathbb{O}$ and $\mathrm{Im}(R)=Fve_{22}+Fve_{12}$. Then, up to conjugation by automorphism, antiautomorphism and up to multiplication by a~scalar, an operator $R$ acts in one of the following ways for some $\alpha\in F$ (an operator $R$ is zero on unspecified basic elements $e_{ij}$, $ve_{ij}$):
\begin{enumerate}
    \item $R(ve_{11}) = ve_{22}$, $R(ve_{21}) = ve_{22} + \alpha ve_{12}$, $\alpha\neq 0$,
    \item $R(e_{21}) = ve_{12}$, $R(ve_{21}) = ve_{22},$
    \item $R(e_{11}) = R(e_{12}) = -R(e_{21}) = - R(e_{22}) = ve_{22} + ve_{12}, R(ve_{11}) = - ve_{12}, R(ve_{21}) = ve_{12}$,
    \item $R(ve_{11}) = ve_{12}$, $R(ve_{21}) = ve_{22} + ve_{12}$,
    \item $R(ve_{11}) =  ve_{12}$, $R(ve_{21}) = \alpha ve_{22}$, $\alpha\neq 0$,
    \item $R(ve_{11}) = ve_{22}$, $R(ve_{21}) = \alpha ve_{12}$, $\alpha\neq 0$.
\end{enumerate}
}

\textbf{Proof.} Note that $\mathrm{Im}(R)$ is a~trivial algebra. Thus, the following equalities hold for any $z\in\mathbb{O}$:
\[ve_{12}R(z)=R(z)ve_{12}=ve_{22}R(z)=R(z)ve_{22}=0.\]
In addition, for any $x,y\in\mathbb{O}$ we have
\[0=R(R(x)y+xR(y)).\]

Let $R(ve_{12})=\gamma_1 ve_{22} + \delta_1 ve_{12}$. Then
\[0=R(ve_{12})R(e_{11}) = R(R(ve_{12})e_{11}) = \delta_1 R(ve_{12}).\]
If $R(ve_{12})=0$, then $\delta_1=0$, so in any case $\delta_1=0$. Further,
\[0=R(ve_{22})R(e_{12})=R(\gamma_1 ve_{22}\cdot e_{12}) = \gamma_1 R(ve_{12}) = \gamma_1^2 ve_{22},\]
whence $\gamma_1 = 0$. Therefore, $R(ve_{12})=0$.

Let $R(ve_{22})=\gamma_2 ve_{22} + \delta_2 ve_{12}$. Then
\[0=R(e_{11})R(ve_{22}) = R(e_{11}R(ve_{22})) = \gamma_2 R(ve_{22}).\]
As above, we have $\gamma_2 = 0$. Further,
\[0 = R(e_{21})R(ve_{12}) = R(e_{21}R(ve_{12})) = - \delta_2 R(ve_{22}) = -\delta_2^2 e_{12},\]
hence $\delta_2 = 0$. Therefore, $R(ve_{22})=0$ and $R^2=0$. 

Let $R(1)=\alpha ve_{22} + \beta ve_{12}$. Then for any $x\in\mathbb{O}$ we have
\[0=R(R(1)x+1R(x))=R(R(1)x).\]
Similarly, $0=R(yR(1))$.
In particular, for $x=ve_{11}$ we obtain
\[0=R(R(1)ve_{11})=R(\alpha e_{11} - \beta e_{12}).\]
For $y=ve_{11}$ we obtain
\[0=R(ve_{11}R(1))=R(\alpha e_{22} + \beta e_{12}).\]
Adding the last two equalities, we obtain $0=\alpha R(1)$.
If $R(1)=0$, then $\alpha = 0$. Thus, in any case $\alpha = 0$. For $x=y=ve_{21}$ we have
\begin{gather*}0=R(ve_{21}R(1))=-\beta R(e_{11}),\\
0=R(R(1)ve_{21})=-\beta R(e_{22}).\end{gather*}
Adding the last two equalities, we obtain $0=-\beta R(1)$. If $R(1)=0$, then $\beta = 0$. Thus, in any case $\beta = 0$ and $R(1)=0$.

Let us define $\alpha_1,\alpha_2,\alpha_3,\beta_1,\beta_2,\beta_3$ as follows.
Let 
\begin{gather}
    \label{e1}R(e_{11}) = \alpha_1 ve_{22} + \beta_1 ve_{12},\\
    \label{e2}R(e_{12}) = \alpha_2 ve_{22} + \beta_2 ve_{12},\\
    \label{e3}R(e_{21}) = \alpha_3 ve_{22} + \beta_3 ve_{12}.
\end{gather}
According to the above we have
\[R(e_{22})= - \alpha_1 ve_{22} - \beta_1 ve_{12}.\]
If $i\neq j$ then
\begin{gather*}
    e_{ii}\cdot ve_{k2} = \delta_{k,3-i}ve_{3-i,2}\in\mathrm{Im}(R),\\
    e_{ij}\cdot ve_{k2} = -\delta_{j,k}ve_{i2}\in\mathrm{Im}(R),\\
    ve_{k2}\cdot e_{ii} = \delta_{i,k}ve_{i2}\in\mathrm{Im}(R),\\
    ve_{k2}\cdot e_{ij} = \delta{j,k}ve_{i2}\in\mathrm{Im}(R).
\end{gather*}
Thus, for any $i,j\in\{1,2\}$ we obtain
\[e_{ij}\mathrm{Im}(R)+\mathrm{Im}(R)e_{ij}\subseteq\mathrm{Im}(R)\subseteq\ker (R).\]

Given the last embedding, the condition that $R$ is a~Rota-Baxter operator is equivalent to the following equalities:
\begin{gather*}
    0=R(R(e_{ij})y) = R(\alpha_m ve_{22} \cdot y + \beta_m ve_{12} \cdot y),\\
    0=R(yR(e_{ij})) = R(\alpha_m y \cdot ve_{22} + \beta_m y\cdot ve_{12}),\\
    0=R(ve_{kl})R(ve_{ij}) = R(R(ve_{ij})ve_{kl}+ve_{ij}R(ve_{kl})),
\end{gather*}
and it is sufficient to consider $y\in vM_2(F)$. Considering the first equality with $y=ve_{kl}$ for $k,l\in\{1,2\}$, we get
\begin{gather*}
    0=R((\alpha_m ve_{22} + \beta_m ve_{12})ve_{kl}) = \alpha_m R(\delta_{l,1}e_{k1})-\beta_m R(\delta_{l,1}e_{k2}).
\end{gather*}
If $l=2$ then we have $0=0$. If $l=1$ and $k=1,2$ then
\begin{gather*}
    0=(\alpha_m\alpha_1-\beta_m\alpha_2)ve_{22}+(\alpha_m\beta_1-\beta_m\beta_2)ve_{12},\\
    0=(\alpha_m\alpha_3+\beta_m\alpha_1)ve_{22}+(\alpha_m\beta_3+\beta_m\beta_1)ve_{12}.
\end{gather*}
Similarly the second equality imply
\begin{gather*}
    0=(-\alpha_m\alpha_1+\beta_m\alpha_2)ve_{22}+(-\alpha_m\beta_1+\beta_m\beta_2)ve_{12}=0,\\
    0=(\alpha_m\alpha_3+\beta_m\alpha_1)ve_{22}+(\alpha_m\beta_3+\beta_m\beta_1)ve_{12}=0.
\end{gather*}
Thus, the condition that $R$ is a~Rota-Baxter operator is equivalent to the following conditions
\begin{gather}
    \label{T1}\alpha_m\alpha_1-\beta_m\alpha_2 = 0,\\
    \label{T2}\alpha_m\beta_1-\beta_m\beta_2 = 0,\\
    \label{T3}\alpha_m\alpha_3+\beta_m\alpha_1 = 0,\\
    \label{T4}\alpha_m\beta_3+\beta_m\beta_1 = 0,\\
    \label{T5}R(R(ve_{ij})ve_{kl}+ve_{ij}R(ve_{kl})) = 0.
\end{gather}
Moreover, in the last equality it is sufficient to consider the pairs $(i,j)$ and $(k,l)$ from the set $\{(1,1),(2,1)\}$, since otherwise the equality $0= 0$.

Let us define $\gamma,\delta,\mu,\nu$ as follows. Let
\begin{gather*}
    R(ve_{11}) = \gamma ve_{22} + \delta ve_{12},\\
    R(ve_{21}) = \mu ve_{22} + \nu ve_{12}.
\end{gather*}

1) Let $\alpha_1=0$. Then the equality \eqref{T3} for $m=3$ implies $\alpha_3=0$, and the equality \eqref{T4} for $m=1$ implies $\beta_1=0$. The last together with the equality \eqref{T2} for $m=2$ implies $\beta_2=0$. The equality \eqref{T4} for $m=2$ implies $\alpha_2\beta_3=0$.

Then we have 
\[0 = R(R(ve_{11})ve_{21}+ve_{11}R(ve_{21})) = \nu\alpha_2 ve_{22} + \gamma\beta_3 ve_{12},\]
whence we get $\nu\alpha_2=\gamma\beta_3=0$.

1.a) Let $\alpha_2=0$. 

1.a.a) Let $\beta_3\neq 0$. Then $\gamma = 0$ and after multiplication by $\frac{1}{\beta_3}$ we obtain
\[R(e_{21}) = ve_{12}, \quad R(ve_{11}) = \delta ve_{12}, \quad R(ve_{21}) = \mu ve_{22} + \nu ve_{12}, \quad \mu\neq 0.\]
 $R$ acts in a~zero manner on the remaining basic elements.

1.a.a.a) Let $\delta\neq 0$. Then Proposition~8 for $\alpha=\delta$ allows us to assume that
\[R(e_{21}) = \frac{1}{\delta} ve_{12},\quad R(ve_{11}) = \frac{1}{\delta} ve_{12},\quad \delta^2 R(ve_{21}) = \mu ve_{22} + \nu ve_{12}, \quad\mu\neq 0.\]
After multiplying by the scalar $\delta$, we get
\[ R(e_{21}) = ve_{12},\quad R(ve_{11}) = ve_{12},\quad R(ve_{21}) = \mu_1 ve_{22} + \nu_1 ve_{12}, \quad \mu_1\neq 0,\]
where $\nu_1 = \frac{\nu}{\delta}$, $\mu_1=\frac{\mu}{\delta}$

\medskip 1.a.a.a.a) Let $\nu\neq 0$. By Proposition~6 for $\alpha = \frac{1}{\nu_1}$ we can assume that
\[R(e_{21}) = \frac{1}{\nu_1^2}ve_{12},\quad R(ve_{11}) = \frac{1}{\nu_1}ve_{12},\quad \nu_1R(ve_{21}) = \mu_1 ve_{22} + ve_{12},\quad \mu_1\neq 0.\]
By Proposition~8 for $\alpha = \nu_1$ we can assume that
\[R(e_{21}) = \frac{1}{\nu_1^3}ve_{12},\quad R(ve_{11}) = \frac{1}{\nu_1^3}ve_{12},\quad \nu_1^3R(ve_{21}) = \mu_1 ve_{22} + ve_{12},\quad \mu_1\neq 0.\]
After multiplication by the scalar $\nu_1^3$ we obtain 
\[R(e_{21}) = ve_{12},\quad R(ve_{11}) = ve_{12},\quad R(ve_{21}) = \mu_1 ve_{22} + ve_{12},\quad \mu_1\neq 0.\]
Conjugation by the automorphism from Proposition~9 with the scalar $\alpha = 1$ allows us to assume that
\[R(e_{21}) = ve_{12},\quad R(ve_{11}) = ve_{12},\quad R(ve_{21}) = \mu_1 ve_{22},\quad \mu_1\neq 0.\]
Proposition~13 allows us to assume that 
\[ R(e_{21}) = ve_{12},\quad R(ve_{11}) = 0,\quad R(ve_{21}) = \mu_1 ve_{22}, \quad \mu_1\neq 0.\]
Proposition~8 with a~scalar $\mu_1$ and a~multiplication by a~scalar $\mu_1$ allow us to assume that
\[R(e_{21}) = ve_{12}, \quad R(ve_{11}) = 0,\quad R(ve_{21}) = ve_{22}.\]
We have obtained case (2) from the statement of the lemma.

\medskip 1.a.a.a.b) Let $\nu = 0$. Then we have
\[R(e_{21}) = ve_{12},\quad R(ve_{11}) = ve_{12},\quad R(ve_{21}) = \mu_1 ve_{22},\quad \mu_1\neq 0.\]
As above, we have obtained case (2) from the statement of the lemma again. 

\medskip 1.a.a.b) Let $\delta = 0$. Then we have
\[R(e_{21}) = ve_{12},\quad R(ve_{11}) = 0,\quad R(ve_{21}) = \mu ve_{22} + \nu ve_{12}, \quad \mu\neq 0.\]
If $\nu\neq 0$ then, as above, we can assume that 
\[R(e_{21}) = ve_{12},\quad R(ve_{21}) = \mu ve_{22} + ve_{12}, \quad \mu\neq 0.\]
Conjugation by the automorphism from Proposition~9 with the scalar $\alpha = 1$ allows us to assume that
\[R(e_{21}) = ve_{12}, \quad R(ve_{21}) = \mu ve_{22}, \mu\neq 0.\]
Conjugation by the automorphism from Proposition~6 with the scalar $\alpha = \mu$ allows us to assume that
\[R(e_{21}) = \mu^2 ve_{12},\quad R(ve_{21}) = \mu^2 ve_{22}, \mu\neq 0.\]
After multiplying by the scalar $\frac{1}{\mu^2}$ we have obtained case (2) from the statement of the lemma again.

\medskip 1.a.b) Let $\beta_3 = 0$. Then $R(M_2(F)) = 0$.

\medskip 1.a.b.a) Let $\gamma = 0$. After multiplying by the scalar $\frac{1}{\delta}$ we have
\[R(ve_{11}) =  ve_{12},\quad R(ve_{21}) = \mu ve_{22} + \nu ve_{12}.\]

\medskip 1.a.b.a.a) Let $\nu\neq 0$. Conjugation by the automorphism from Proposition~7 with the scalar $\alpha = \frac{1}{\nu}$ and multiplication by the scalar $\frac{1}{\nu}$ allow us to assume that
\[R(ve_{11}) = ve_{12},\quad R(ve_{21}) = \mu ve_{22} + ve_{12}.\]
If $\mu = 1$ then we have obtained case (4) from the statement of the lemma. If $\mu\neq 1$ then Proposition~16 with a~scalar $\frac{1}{\alpha-1}$ gives us case (5) from the statement of the lemma.

\medskip 1.a.b.a.b) Let $\nu = 0$. Then
\[ R(ve_{11}) =  ve_{12},\quad R(ve_{21}) = \mu ve_{22}.\]
We have obtained case (5) from the statement of the lemma.

\medskip 1.a.b.b) Let $\gamma\neq 0$. Then we can assume that
\[R(ve_{11}) = ve_{22} + \delta ve_{12},\quad R(ve_{21}) = \mu ve_{22} + \nu ve_{12}.\]

\medskip 1.a.b.b.a) Let $\delta = 0$ and $\mu\neq 0$. Conjugation by the automorphism from Proposition~7 with the scalar $\alpha = \frac{1}{\mu}$ and multiplication by the scalar $\frac{1}{\mu}$ allow us to assume that
\[R(ve_{11}) = ve_{22},\quad R(ve_{21}) = ve_{22} + \nu_1 ve_{12}.\]
We have obtained case (1) from the statement of the lemma.

\medskip 1.a.b.b.b) Let $\delta=\mu=0.$ Then
\[R(ve_{11}) = ve_{22},\quad R(ve_{21}) = \nu ve_{12}.\]
We have obtained case (6) from the statement of the lemma.

\medskip 1.a.b.b.c) Let $\delta\neq 0$ and $\mu\neq 0$. Conjugation by the automorphism from Proposition~6 with the scalar $\alpha = \frac{1}{\delta}$ allows us to assume that
\[ R(ve_{11}) = ve_{22} + ve_{12},\quad R(ve_{21}) = \mu_1 ve_{22} + \nu_1 ve_{12}.\]
Since $\mu_1\neq \nu_1$, then Proposition~16 with a~scalar $\alpha = 1$ allows us to assume that
\[ R(ve_{11}) = ve_{22},\quad R(ve_{21}) = \mu_1' ve_{22} + \nu_1' ve_{12},\]
where $\nu_1'\neq 0$. A conjugation by an automorphism from Proposition~7 with a~scalar $\frac{1}{\mu_1'}$ and a~multiplication by a~scalar $\frac{1}{\mu_1'^2}$ give us case (1) from the statement of the lemma.

\medskip 1.b) Let $\alpha_2\neq 0$, then $\beta_3=\nu=0$. By Proposition~13 this case is antiisomorphic to the case 1.a). 

\medskip 2) Let $\alpha_1\neq 0$. Then the equality \eqref{T3} for $m=1$ implies $\alpha_3=-\beta_1$. The equality \eqref{T1} for $m=1$ implies $\alpha_1^2=\beta_1\alpha_2$. We have $\alpha_2\neq 0$, because otherwise $\alpha_1=0$. Then the equality \eqref{T1} for $m=2$ implies $\alpha_1=\beta_2$. Particularly, the equalities \eqref{T2} and \eqref{T3} are equivalent. The equality \eqref{T1} for $m=3$ implies $\alpha_2\beta_3=\alpha_1\alpha_3$. The equality \eqref{T4} for $m=1$ implies $\beta_1^2+\alpha_1\beta_3=0$. We have $\alpha_3\neq 0$, because otherwise $\beta_1 = -\alpha_3 = 0$ and $\alpha_1^2 =0$. So, we have
\[\beta_1=-\alpha_3,\quad  \beta_2=\alpha_1,\quad
    \beta_3=-\frac{\alpha_3^2}{\alpha_1},\quad
    \alpha_2 = -\frac{\alpha_1^2}{\alpha_3}.\]
It is easy to see that these conditions are sufficient to satisfy \eqref{T1}--\eqref{T4}. Since we describe operators up to multiplication by a~scalar, we can assume that $\alpha_1=1$. Let us denote $t=\alpha_3$. Then the equalities \eqref{e1}--\eqref{e3} are converted to
\[
    R(e_{11}) = ve_{22} - t ve_{12},\quad R(e_{12}) = -\frac{1}{t} ve_{22} + ve_{12},\quad R(e_{21}) = t ve_{22} - t^2 ve_{12}.
\]
Since we are describing operators up to multiplication by a~scalar, we can assume that we are considering the operator $tR$. Then
\begin{gather}
   R(e_{11}) = t ve_{22} - t^2 ve_{12},\\
   R(e_{12}) = -ve_{22} + t ve_{12},\\
   R(e_{21}) = t^2 ve_{22} - t^3 ve_{12}.
\end{gather}
Further, it is easy to see that
\begin{gather*}
    R(R(ve_{11})ve_{11}+ve_{11}R(ve_{11}))=\gamma R(1) = 0,\\
    R(R(ve_{21})ve_{21}+ve_{21}R(ve_{21}))=-\nu R(1) = 0.
\end{gather*}
Let us consider the equality \eqref{T5} with $(i,j)=(1,1)$ and $(k,l)=(2,1)$. Then we have
\begin{multline*}
R(R(ve_{11})ve_{21}+ve_{11}R(ve_{21})) = R(\gamma e_{21} + \nu e_{12} + (\mu-\delta)e_{22}) = \\
= \gamma (t^2 ve_{22} - t^3 ve_{12}) + \nu (-ve_{22} + tve_{12}) +(\mu-\delta) (-tve_{22}+t^2 ve_{12}),
\end{multline*}
which is equivalent to the condition 
\[\gamma t^2 - \nu + (\delta - \mu) t = 0.\]
The parameters $(i,j)=(2,1)$ and $(k,l)=(1,1)$ imply the same condition. It means that conditions \eqref{T1}--\eqref{T5} are equivalent to the following equalities:
\begin{gather*}
  R(e_{11}) = t ve_{22} - t^2 ve_{12}, \quad  R(e_{12}) = -ve_{22} + t ve_{12},\quad R(e_{21}) = t^2 ve_{22} - t^3 ve_{12},\\
  R(e_{22}) = -t ve_{22} + t^2 ve_{12}, \quad R(ve_{11}) = \gamma ve_{22} + \delta ve_{12},\\ R(ve_{21}) = \mu ve_{22} + (\gamma t^2 + (\delta-\mu)t)ve_{12}.
\end{gather*}
After conjugation by the automorphism from Proposition~6 with $\alpha = -\frac{1}{t}$ and multiplication by the scalar $\frac{1}{t}$ we can assume that
\begin{gather*}
  R(e_{11}) = ve_{22} + ve_{12},\quad R(e_{12}) = ve_{22} + ve_{12},\quad R(e_{21}) = -ve_{22} - ve_{12},\\
  R(e_{22}) = -ve_{22} - ve_{12}, \quad R(ve_{11}) = \gamma_1 ve_{22} - \delta_1 ve_{12},\\
  R(ve_{21}) = \mu_1 ve_{22} + (\gamma_1 + \delta_1 + \mu_1) ve_{12},
\end{gather*}
where $\gamma_1 = \frac{\gamma}{t}$, $\delta_1 = \frac{\delta}{t^2}$, $\mu_1 = -\frac{\mu}{t^2}$

2.a) If $\gamma_1 = 0$, then a~conjugation by the automorphism from Proposition~8 with $\alpha = \delta_1$ and multiplication by the scalar $\delta_1$ allow us to assume that
\begin{gather*}
  R(e_{11}) = ve_{22} + ve_{12},\quad R(e_{12}) = ve_{22} + ve_{12},\quad R(e_{21}) = -ve_{22} - ve_{12},\\
  R(e_{22}) = -ve_{22} - ve_{12}, \quad R(ve_{11}) = -ve_{12},\quad R(ve_{21}) = \mu_2 ve_{22} + (1 + \mu_2) ve_{12},
\end{gather*}
where $\mu_2 = \frac{\mu_1}{\delta_1}$. Conjugation by the automorphism from Proposition~10 with $\alpha = -\mu_2$ gives us the case (3) from the statement of the lemma.

2.b) If $\gamma_1\neq 0$, then conjugation by the automorphism from Proposition~8 with $\alpha = \gamma_1$ and multiplication by the scalar $\gamma_1$ allow us to assume that
\begin{gather*}
  R(e_{11}) = ve_{22} + ve_{12},\quad R(e_{12}) = ve_{22} + ve_{12},\quad R(e_{21}) = -ve_{22} - ve_{12},\\
  R(e_{22}) = -ve_{22} - ve_{12}, \quad R(ve_{11}) = ve_{22} - \delta_2 ve_{12},\\ R(ve_{21}) = \mu_2 ve_{22} + (1 + \delta_2 + \mu_2) ve_{12},
\end{gather*}
where $\mu_2 = \frac{\mu_1}{\gamma_1}$, $\delta_2 = \frac{\delta_1}{\gamma_1}$. We have $\delta_2\neq -1$ because $\dim(\mathrm{Im}(R))=2$. Conjugation by the automorphism from Proposition~9 with $\alpha = -\frac{1}{2}$ give us
\begin{gather*}
  R(e_{11}) = ve_{22} + ve_{12},\quad   R(e_{12}) = ve_{22} + ve_{12},\quad R(e_{21}) = -ve_{22} - ve_{12},\\
  R(e_{22}) = -ve_{22} - ve_{12}, \quad  R(ve_{11}) = -(\delta_2 + 1) ve_{12},\\
  R(ve_{21}) = (\mu_2 + 1/2) ve_{22} + (1 + \delta_2 + \mu_2 + 1/2) ve_{12}.
\end{gather*}
We have the case $\gamma_1=0$ which reduces to case (3) from the statement of the lemma. The lemma is proven.

\medskip\textbf{Corollary 2.} \textit{Let $R$ be a~Rota-Baxter operator of zero weight on the split Cayley-Dickson algebra $\mathbb{O}$ and $\mathrm{Im}(R)=Fve_{22}+Fve_{12}$. If a~field $F$ is quadratically closed then, up to conjugation by automorphism, antiautomorphism and up to multiplication by a~scalar, an operator $R$ acts in one of the following ways for some $\alpha\in F$ (an operator $R$ is zero on unspecified basic elements $e_{ij}$, $ve_{ij}$):
\begin{enumerate}
    \item $R(e_{21}) = ve_{12}$, $R(ve_{21}) = ve_{22},$
    \item $R(e_{11}) = R(e_{12}) = -R(e_{21}) = - R(e_{22}) = ve_{22} + ve_{12}, R(ve_{11}) = - ve_{12}, R(ve_{21}) = ve_{12}$,
    \item $R(ve_{11}) = ve_{12}$, $R(ve_{21}) = ve_{22} + ve_{12}$,
    \item $R(ve_{11}) =  ve_{12}$, $R(ve_{21}) = \alpha ve_{22}$, $\alpha\neq 0$,
\end{enumerate}
}

\textbf{Proof.} Let us consider the cases from the statement of Lemma~4.

(1) Let $\beta$ be a~root of the equation $-\beta^2\alpha - \beta + 1 = 0$. Then Proposition~17 with a~scalar $\beta$ allows us to assume that
\[R(ve_{11}) = -\alpha\beta ve_{12}, \quad R(ve_{21}) = (1+\alpha\beta) ve_{22} + \alpha ve_{12}.\]
By Proposition~7 with a~scalar $-\alpha\beta$ we can assume that $R(ve_{11}) = ve_{12}$, $R(ve_{21}) = \mu ve_{22} + \nu ve_{12}$ for some $\mu,\nu\in F$. We have a~case 1.a.b.a) from the proof of Lemma~4, which was proved to be equivalent to the cases (4) and (5) from the statement of Lemma~4 (cases (3) and (4) from the current corollary).

(6) This case is fixed by Proposition~6 with the scalar $\frac{1}{\sqrt{\alpha}}$, we have obtained $R(ve_{11})=ve_{22}$, $R(ve_{21})=ve_{12}$. Proposition~17 with a~scalar $\alpha = -1$ allows us to assume that $R(ve_{11})=ve_{12}$, $R(ve_{21}) = ve_{12}-ve_{22}$. Proposition~16 with a~scalar $-\frac{1}{2}$ allows us to assume that $R(ve_{11}) = ve_{12}$, $R(ve_{21}) = -ve_{22}$. We have case (4) from the statement of the current corollary. The corollary is proven.

\section{RB-Operators with three-dimensional image}

In \cite{Octonions} it was proved that there are only two three-dimensional non-unital subalgebras $B$ in $\mathbb{O}$, up to action of automorphism: nilpotent $Fe_{12}+Fve_{12}+Fve_{22}$ and idempotent $Fe_{11}+Fve_{12}+Fve_{22}$. Let us describe the Rota-Baxter operators of zero weight on $\mathbb{O}$ with these images.

\medskip \textbf{Lemma 5.} \textit{Let $R$ be a~Rota-Baxter operator on the split Cayley-Dickson algebra~$\mathbb{O}$ and $\mathrm{Im}(R)=Fe_{12} + Fve_{12} + Fve_{22}$. Then \[\ker (R) = L(e_{11},e_{12},e_{22},ve_{12},ve_{22}),\]
and, up to conjugation by automorphism, antiautomorphism and up to multiplication by a~scalar, an operator $R$ acts in one of the following ways for some $\alpha\in F$:
\begin{enumerate}
\item $R(e_{21})=\alpha e_{12},\quad R(ve_{11})=ve_{12},\quad R(ve_{21})=ve_{22} + ve_{12}$,
where $\alpha\neq 0$;
\item $R(e_{21}) = \alpha e_{12},\quad R(ve_{11})=ve_{12},\quad R(ve_{21})=ve_{22},$
where $\alpha\neq 0$;
\item $R(e_{21})= e_{12},\quad R(ve_{11})=ve_{12},\quad R(ve_{21})=ve_{22} + e_{12};$
\item $R(e_{21})=\alpha e_{12},\quad R(ve_{11})=ve_{12},\quad R(ve_{21})=ve_{22} + ve_{12} + e_{12},$
where $\alpha\neq 0$.
\end{enumerate}
}

\textbf{Proof.} Let us notice that $(\mathrm{Im}(R))^2=L(ve_{12})$ and we have the following equalities on $\mathbb{O}$: 
\[ve_{12}R(z)=R(z)ve_{12}=0.\]
Particularly there is an embedding for any $x,y\in\mathbb{O}$:
\[R(R(x)y+xR(y))\in L(ve_{12}).\]
Let $R(ve_{12})=\xi ve_{22} + \eta ve_{12} + \mu e_{12}$. Then we have
\[L(ve_{12})\ni R(ve_{12})R(e_{12}) = R(R(ve_{12})e_{12}) = \xi R(ve_{12}) = \xi^2 ve_{22} + \xi\eta ve_{12} + \xi\mu e_{12},\]
whence $\xi = 0$. Further,
\[L(ve_{12})\ni R(ve_{12})R(ve_{22})=R(R(ve_{12})ve_{22}) = R(-\mu ve_{12}) = -\mu\eta ve_{12} - \mu^2 e_{12},\]
whence $\mu = 0$. So, $R(ve_{12})=\eta ve_{12}$. Particularly, $ R(z)R(ve_{12})=0$ for any $z\in\mathbb{O}$. Then we have
\[0=R(e_{22})R(ve_{12})=R(\eta e_{22}\cdot ve_{12}) = \eta R(ve_{12}) = \eta^2 ve_{12},\]
whence $\eta=0$. So, $R(ve_{12})=0$.

1) Let $R^2=0$. Then $R(ve_{22})=R(e_{12})=0$. Let $R(1)=\alpha ve_{22} + \beta ve_{12} + \gamma e_{12}$. Then $R(R(1)x)=R(1)R(x)\in L(ve_{12})$ and $R(yR(1))\in L(ve_{12})$. For $x=ve_{11}$ we obtain
\[L(ve_{12})\ni R(R(1)ve_{11})= R(\alpha e_{11} - \beta e_{12}).\]
For $y=ve_{11}$ we obtain
\[L(ve_{12})\ni R(ve_{11}R(1)) = R(\alpha e_{22} + \beta e_{12}).\]
Adding these equalities, we get $\alpha R(1) \in L(ve_{12})$, i.e. $\alpha = 0$. 

If $x=y=ve_{21}$ then we have
\begin{gather*}
L(ve_{12})\ni R(R(1)ve_{21}) = -\beta R(e_{22}) - \gamma R(ve_{11}),\\
L(ve_{12})\ni R(ve_{21}R(1)) = -\beta R(e_{11}) + \gamma R(ve_{11}).
\end{gather*}
Adding these equalities, we get $-\beta R(1) \in L(ve_{12})$, whence $\beta\gamma = 0$. If $\beta = 0$, then $R(1)=\gamma e_{12}$ and
\[L(ve_{12})\ni R(e_{11})R(1) = R(e_{11}R(1)) = \gamma R(e_{12}) = \gamma^2 e_{12},\]
whence $\gamma = 0$. Thus in any case $\gamma = 0$ and $R(1) = \beta ve_{12}$. Besides, it means that $R(1)R(x)=0$ for any $x\in\mathbb{O}$. Then from the equalities above we have $\beta R(e_{11}), \beta R(e_{22}) = 0$. So, $\beta R(1) = 0 $ and $\beta^2 = 0$.

We have
\[ve_{22},ve_{12},e_{12},1\in\ker (R).\]
Let
\begin{gather*}
R(e_{11})=\alpha_1 ve_{22} + \beta_1 ve_{12} + \gamma_1 e_{12},\\
R(e_{21})=\alpha_2 ve_{22} + \beta_2 ve_{12} + \gamma_2 e_{12},\\
R(ve_{11})=\alpha_3 ve_{22} + \beta_3 ve_{12} + \gamma_3 e_{12},\\
R(ve_{21})=\alpha_4 ve_{22} + \beta_4 ve_{12} + \gamma_4 e_{12}.
\end{gather*}
Then
\[L(ve_{12})\ni R(ve_{11})R(e_{11}) = \alpha_1 R(e_{22}) = -\alpha_1^2 ve_{22} -\alpha_1\beta_1ve_{12}-\alpha_1\gamma_1 e_{12},\]
whence $\alpha_1^2 = 0$ and $\alpha_1=0$. Further,
\[L(ve_{12})\ni R(e_{11})R(e_{21}) = \gamma_1R(e_{11}) = \beta_1\gamma_1 ve_{12} + \gamma_1^2 e_{12},\]
whence $\gamma_1^2 = 0$ and $\gamma_1=0$. So, $R(e_{11})R(x)=0$ for any $x\in\mathbb{O}$, because $R(e_{11})=\beta_1 ve_{12}$. Then
\[0=R(e_{11})R(ve_{21})=R(-\beta_1e_{22}+\alpha_4 ve_{22}+\gamma_4 e_{12}) = -\beta_1 R(e_{22}) = \beta_1^2 ve_{12},\]
whence $\beta_1=0$. Thus, $R(e_{11})=0$ and we found the kernel of operator $R$:
\[\ker (R)=L(e_{11},e_{12},e_{22},ve_{12},ve_{22}).\]

It is easy to see that for this $\ker (R)$ the condition $R(x)R(y)=R(R(x)y+xR(y))$ is executed automatically, if $x\in\ker (R)$ or $y\in\ker (R)$. Indeed, let $x\in \ker (R)$. Then this condition is equivalent to $0=R(xR(y))$. The last condition is executed because $\mathrm{Im}(R)\subset \ker (R)$ and $\ker (R)$ is a~subalgebra. The case $y\in\ker (R)$ is similar.

Thus, for $R$ to be a~Rota-Baxter operator, it is necessary and sufficient that $R(x)R(y)=R(R(x)y+xR(y))$ holds for $x, y\in\{e_{21},ve_{11},ve_{21}\}$. Let us notice that $e_{21}\mathbb{O}=L(e_{21},e_{22},ve_{21},ve_{22})\subset \ker (R)$, whence $R(e_{21}y)=0$ for any $y\in\mathbb{O}$. Then
\[
(\alpha_2\gamma_4-\alpha_4\gamma_2)ve_{12}=R(e_{21})R(ve_{21}) = (\alpha_2^2-\gamma_2\alpha_3)ve_{22}+(\alpha_2\beta_2-\gamma_2\beta_3)ve_{12}+\gamma_2(\alpha_2-\gamma_3)e_{12},\]
whence
\begin{gather}
\label{5.1} \alpha_2^2-\gamma_2\alpha_3 = \alpha_2\gamma_2-\gamma_2\gamma_3 = 0,\\
\label{5.2} \alpha_2(\gamma_4-\beta_2)=\gamma_2(\alpha_4-\beta_3).
\end{gather}
Further,
\[(\alpha_2\gamma_3-\alpha_3\gamma_2)ve_{12}=R(e_{21})R(ve_{11})=R(R(e_{21})ve_{11})=R(\alpha_2 e_{11} - \beta_2 e_{12}) = 0,\]
whence 
\begin{equation}\label{5.3}
\alpha_2\gamma_3-\alpha_3\gamma_2=0.
\end{equation}
Finally,
\[
(\alpha_3\gamma_4-\alpha_4\gamma_3)ve_{12} = R(ve_{11})R(ve_{21}) =
(\alpha_3\alpha_2-\gamma_3\alpha_3)ve_{22} + (\alpha_3\beta_2-\gamma_3\beta_3)ve_{12} + (\alpha_3\gamma_2-\gamma_3\gamma_3)e_{12},
\]
whence 
\begin{gather} \label{5.4} \alpha_3(\alpha_2-\gamma_3)=0, \quad \alpha_3\gamma_2-\gamma_3^2 = 0,\\
\label{5.5} \alpha_3(\gamma_4-\beta_2) = \gamma_3(\alpha_4-\beta_3).
\end{gather}

1.a) Let us assume that $\alpha_3=0$. Then by \eqref{5.4} $\alpha_2=\gamma_3=0$ and by \eqref{5.2} $\gamma_2(\alpha_4-\beta_3)=0$. Since $\gamma_2\neq 0$ (otherwise $\dim\mathrm{Im}(R) < 3$), then $\alpha_4=\beta_3\neq 0$. Then, up to multiplication by a~scalar, an operator $R$ acts as follows:
\begin{gather*}
R(e_{21})=\beta_2 ve_{12} + \gamma_2 e_{12},\quad R(ve_{11})=ve_{12},\\
R(ve_{21})=ve_{22} + \beta_4 ve_{12} + \gamma_4 e_{12},
\end{gather*}
where $\gamma_2\neq 0$.

A conjugation by automorphism from Proposition~9 with a~scalar $\alpha = -\frac{\beta_2}{\gamma_2}$ allows us to assume that
\begin{gather*}
R(e_{21})=\gamma_2 e_{12},\quad R(ve_{11})=ve_{12},\\
R(ve_{21})=ve_{22} + \beta' ve_{12} + \gamma' e_{12}.
\end{gather*}

1.a.a) Let us assume that $\gamma' = 0$ and $\beta'\neq 0$. After a~conjugation by an automorphism from Proposition~7 with a~scalar $\frac{1}{\beta'}$ and after a~multiplication by $\frac{1}{\beta'}$ we obtain
\[ R(e_{21})=\gamma'' e_{12},\quad R(ve_{11})=ve_{12},\quad R(ve_{21})=ve_{22} + ve_{12}.\]

1.a.b) Let us assume that $\gamma' = \beta' = 0$. Then
\[R(e_{21})=\gamma_2 e_{12},\quad R(ve_{11})=ve_{12},\quad R(ve_{21})=ve_{22}.\]

1.a.c) Let us assume that $\gamma'\neq 0$. Then after a~conjugation by an automorphism from Proposition~8 with a~scalar $\frac{1}{\gamma'}$ and after a~multiplication by a~scalar $\frac{1}{\gamma'^2}$ we obtain
\[R(e_{21})=\gamma_2' e_{12},\quad R(ve_{11})=ve_{12},\quad R(ve_{21})=ve_{22} + \beta' ve_{12} + e_{12}.\]

1.a.c.a) Let us assume that $\beta' = 0$. Then we have
\[R(e_{21})=\gamma_2' e_{12},\quad R(ve_{11})=ve_{12},\quad R(ve_{21})=ve_{22} + e_{12}.\]
A conjugation by the automorphism from the Proposition~8 with a~scalar $\frac{1}{\gamma_2'^2}$ allows us to assume that
\[R(e_{21})=\gamma_2' e_{12},\quad R(ve_{11})= \gamma_2'^4 ve_{12},\quad R(ve_{21})= \gamma_2'^4 ve_{22} + \gamma_2'^2 e_{12}.\]
A conjugation by the automorphism from the Proposition~7 with a~scalar $\gamma_2'$ and a~multiplication by a~scalar $\frac{1}{\gamma_2'^3}$ allow us to assume that
\[R(e_{21})= e_{12},\quad R(ve_{11})= ve_{12},\quad R(ve_{21})= ve_{22} + e_{12}.\]

1.a.c.b) Let us assume that $\beta'\neq 0$. A conjugation by the automorphism from the Proposition~7 with a~scalar $\frac{1}{\beta'}$ and multiplication by a~scalar $\frac{1}{\beta'}$ allow us to assume that
\[ R(e_{21})=\frac{\gamma_2'}{\beta'^3} e_{12},\quad R(ve_{11})=ve_{12},\quad R(ve_{21})=ve_{22} + ve_{12} + \frac{1}{\beta'^2}e_{12}.\]
As usual, by the Proposition~8 with a~scalar $\beta'^2$ we can assume that
\[R(e_{21})=\gamma_2'\beta' e_{12},\quad R(ve_{11})=ve_{12},\quad R(ve_{21})=ve_{22} + ve_{12} + e_{12}.\]

\medskip 1.b) Let us assume that $\alpha_3\neq 0$. Then by \eqref{5.4} and \eqref{5.5} we have
\[\alpha_2 = \gamma_3, \quad \beta_2 = \frac{\alpha_3\gamma_4-\gamma_3(\alpha_4-\beta_3)}{\alpha_3} ,\quad \gamma_2 = \frac{\gamma_3^2}{\alpha_3}.\]
Up to multiplication by a~scalar, we can assume that $\alpha_3 = 1$ and
\begin{gather*}
R(e_{21})=\gamma_3 ve_{22} + (\gamma_4-\gamma_3(\alpha_4-\beta_3)) ve_{12} + \gamma_3^2 e_{12},\\
R(ve_{11})= ve_{22} + \beta_3 ve_{12} + \gamma_3 e_{12},\\
R(ve_{21})=\alpha_4 ve_{22} + \beta_4 ve_{12} + \gamma_4 e_{12},
\end{gather*}
where $\gamma_4-\gamma_3\alpha_4 \neq 0$.

1.b.a) Let us assume that $\gamma_3 = 0$. Then
\begin{gather*}
R(e_{21})=\gamma_4 ve_{12},\quad R(ve_{11})= ve_{22} + \beta_3 ve_{12},\\
R(ve_{21})=\alpha_4 ve_{22} + \beta_4 ve_{12} + \gamma_4 e_{12},
\end{gather*}
where $\gamma_4\neq 0$. A conjugation by the involution from Proposition~5 allows to assume that
\begin{gather*}
R(e_{21})=\beta_3 ve_{12} + e_{12},\quad R(ve_{11})= \gamma_4 ve_{12},\\
R(ve_{21})=\gamma_4 ve_{22} + \beta_4 ve_{12} + \alpha_4 e_{12},
\end{gather*}
i.e. we have the case (1.a) with $\alpha_3 = 0$. 

1.b.b) Let us assume that $\gamma_3 \neq 0$. Then the conjugation by automorphism from Proposition~9 with a~scalar $\alpha = -\frac{\gamma_4 - \gamma_3 (\alpha_4-\beta_3)}{\gamma_3^2}$ allows us to assume that
\begin{gather*}
R(e_{21})=\gamma_3 ve_{22} + \gamma_3^2 e_{12},\\
R(ve_{11})= ve_{22} + (\alpha_4 - \gamma_4/\gamma_3) ve_{12} + \gamma_3 e_{12},\\
R(ve_{21})=\alpha_4 ve_{22} + \beta_4' ve_{12} + \gamma_4 e_{12}.
\end{gather*}
Conjugation by automorphism from Proposition~6 with a~scalar $\alpha = \frac{1}{\gamma_3}$ allows us to assume that
\begin{gather*}
R(e_{21})=ve_{22} + e_{12},\\
R(ve_{11})= ve_{22} + (\alpha_4' - \gamma_4') ve_{12} + e_{12},\\
R(ve_{21})=\alpha_4' ve_{22} + \beta_4'' ve_{12} + \gamma_4' e_{12}.
\end{gather*}
where $\alpha_4' = \frac{\alpha_4}{\gamma_3}$, $\beta_4''=\frac{\beta_4'}{\gamma_3^2}$, $\gamma_4' = \frac{\gamma_4}{\gamma_3^2}$.
A conjugation by the automorphism from the Proposition~13 with a~scalar $\alpha = 1$ allows us to assume that
\begin{gather*}
R(e_{21})=e_{12},\quad R(ve_{11})= (\alpha_4' - \gamma_4') ve_{12},\\
R(ve_{21})=(\alpha_4'-\gamma_4') ve_{22} + \beta_4'' ve_{12} + \gamma_4' e_{12}.
\end{gather*}
So, we have a~case $\alpha_3 = 0$, it was considered in (1.a)

\medskip 2) Let us denote $V=\mathrm{Im}(R)\cap\ker (R)$. Let $\dim V= 2$. It means that $\dim(\mathrm{Im}(R^2)) = 1$. It has been proved above that $ve_{12}\in V$. Let
\begin{gather*}
R(1)=\alpha ve_{22} + \beta ve_{12} + \gamma e_{12},\\
R(e_{12})=\alpha_1 ve_{22} + \beta_1 ve_{12} + \gamma_1 e_{12},\\
R(ve_{22})=\alpha_2 ve_{22} + \beta_2 ve_{12} + \gamma_2 e_{12},\\
R(ve_{11}) = \alpha_3 ve_{22} + \beta_3 ve_{12} + \gamma_3 e_{12},\\
R(ve_{21}) = \alpha_4 ve_{22} + \beta_4 ve_{12} + \gamma_4 e_{12},\\
R(e_{11}) = \alpha_5 ve_{22} + \beta_5 ve_{12} + \gamma_5 e_{12},\\
R(e_{21}) = \alpha_6 ve_{22} + \beta_6 ve_{12} + \gamma_6 e_{12}.
\end{gather*}
Then
\[0=R(1)R(1) = 2R(R(1)) =  2(\alpha\alpha_1+\gamma\alpha_2)ve_{22} + 2(\alpha\beta_1+\gamma\beta_2)ve_{12} + 2(\alpha\gamma_1+\gamma\gamma_2)e_{12},\]
whence 
\begin{equation}\label{5.6}
    \alpha\alpha_1 + \gamma\alpha_2 = 0, \quad \alpha\beta_1+\gamma\beta_2 = 0, \quad \alpha\gamma_1+\gamma\gamma_2 = 0.
\end{equation}
Besides,
\begin{gather*}
0 = R(e_{12})R(1)  = (\alpha_1\alpha_2+\gamma_1\alpha_1)ve_{22} + (\alpha_1\beta_2+\gamma_1\beta_1)ve_{12} + (\alpha_1\gamma_2+\gamma_1\gamma_1)e_{12},\\
0 = R(ve_{22})R(1)  = (\alpha_2\alpha_2+\gamma_2\alpha_1)ve_{22} + (\alpha_2\beta_2+\gamma_2\beta_1)ve_{12} + (\alpha_1\gamma_2+\gamma_2\gamma_1)e_{12},\\
(\alpha_1\gamma_2-\alpha_2\gamma_1)ve_{12} = R(e_{12})R(ve_{22}) \in L(R(ve_{12})) = 0,
\end{gather*}
whence
\begin{gather*}
\alpha_1\gamma_2=\alpha_2\gamma_1, \quad \alpha_1(\alpha_2+\gamma_1)=0, \quad \alpha_1\beta_2+\gamma_1\beta_1 = 0,\\
\alpha_1\gamma_2 + \gamma_1^2 = 0, \quad \alpha_2^2+\gamma_2\alpha_1 = 0, \quad \alpha_2\beta_2+\gamma_2\beta_1 = 0, \gamma_2(\alpha_1+\gamma_1) = 0.
\end{gather*}
Let $\alpha_1 = 0$. Then $\gamma_1 = 0$, $\alpha_2 = 0$, $\gamma_2\beta_1 = 0$. If $\beta_1=0$, then $e_{12}\in V$ and $R(ve_{22})=\beta_2 ve_{12} + \gamma_2 e_{12}$. If $\beta_1 \neq 0$,  then $\gamma_2 = 0$ and $R(e_{12}) = \beta_1 ve_{12}$, $R(ve_{22}) = \beta_2 ve_{12}$.

If $\alpha_1\neq 0$, then $\alpha_2=-\gamma_1$. If $\gamma_2 = 0$, then $\alpha_2 = -\gamma_1 = 0$, whence $\beta_2 = 0$ and $R(ve_{22})=0$. If $\gamma_2\neq 0$, then $\alpha_1=-\gamma_1 = \alpha_2$. Since $\alpha_1\gamma_2+\gamma_1^2=0$, then $\gamma_1(\gamma_1-\gamma_2)=0$. If $\gamma_1 = 0$, then $\gamma_2\beta_1 = 0$, whence $\beta_1=0$ and $R(e_{12})=0$. If $\gamma_1\neq 0$, then $\gamma_1=\gamma_2=-\alpha_1=-\alpha_2\neq 0$. Thus, $\alpha_1(\beta_2-\beta_1)=0$ and $\beta_1=\beta_2$. Thus, 
\[R(e_{12}) = \alpha_1 ve_{22} + \beta_1 ve_{12} - \alpha_1 e_{12} = R(ve_{22}).\]

Therefore, 4 cases are possible:

2.a) $V=L(e_{12},ve_{12})$, $R(ve_{22}) = \beta_2 ve_{12} + \gamma_2 e_{12}$. It means that $\gamma\gamma_2=0$ and $\gamma\beta_2=0$, whence $\gamma=0$. Then we have
\begin{gather*}
L(ve_{12})\ni R(1)R(ve_{11}) = R(R(1)ve_{11} + R(ve_{11})) = R(\alpha e_{11} + \alpha_3 ve_{22}),\\
L(ve_{12})\ni R(ve_{11})R(1) = R(\alpha_3 ve_{22} + \alpha e_{22}).
\end{gather*}
A sum of these embedding gives us $\alpha R(1) + 2\alpha_3 R(ve_{22}) \in L(ve_{12})$, whence $\alpha^2 = 0$ and $\alpha=0$. Besides, it means that $R(1)R(ve_{11}) = 0$, so $\alpha_3 R(ve_{22}) = 0$, whence $\alpha_3=0$.
Further,
\begin{gather*}
0 = R(1)R(ve_{21}) = R(R(1)ve_{21} + R(ve_{21})) = R(- \beta e_{22} + \alpha_4 ve_{22}),\\
0 = R(ve_{21})R(1) = R(R(ve_{21}) + ve_{21}R(1)) = R(\alpha_4 ve_{22} - \beta e_{11}).
\end{gather*}
Adding these equalities, we get $-\beta R(1) + 2\alpha_4 R(ve_{22}) = 0$, whence $-\beta^2 + 2\alpha_4\beta_2 = 0$. Further,
\begin{gather*}
0 = R(1)R(e_{21}) = R(\beta ve_{22} + R(e_{21})),\\
0 = R(e_{21})R(1) = R(R(e_{21}) - \beta ve_{22}).
\end{gather*}
Subtracting these equalities, we get $2\beta R(ve_{22}) = 0$, whence $\beta\beta_2=0$. So, in any case $\beta = 0$ and $R(1)=0$. 
Since $\ker (R)$ is a~$\mathrm{Im}(R)$-bimodule, then it means that $ve_{22}\in\ker (R)$. It is a~contradiction.

\medskip 2.b) $V=L(ve_{22},ve_{12})$, $R(e_{12}) = \alpha_1 ve_{22} + \beta_1 ve_{12}$. By Proposition~5 this operator is antiisomorphic to the last. So, it is impossible.

\medskip 2.c) $R(e_{12}) = \beta_1 ve_{12}$, $R(ve_{22}) = \beta_2 ve_{12}$, $\beta_1,\beta_2\neq 0$. Further,
\begin{gather*}
L(ve_{12})\ni R(1)R(ve_{11}) = R(R(1)ve_{11} + R(ve_{11})) = R(\alpha e_{11} - \beta e_{12} + \alpha_3 ve_{22} + \gamma_3 e_{12}),\\
L(ve_{12})\ni R(ve_{11})R(1) = R(\alpha_3 ve_{22} + \gamma_3 e_{12} + \alpha e_{22} + \beta e_{12}).
\end{gather*}
Adding these equalities, we get $\alpha R(1) + 2\alpha_3 R(ve_{22}) + 2\gamma_3 R(e_{12})\in L(ve_{12})$, what implies $\alpha^2 = 0$ and $\alpha = 0$. So, by \eqref{5.6} we have $\gamma\beta_2 = 0$, whence $\gamma = 0$. It means that $R(1)R(x) = R(x)R(1) = 0$ for any $x\in\mathbb{O}$. Further, as in the first case, we have
\begin{gather*}
0= R(1)R(e_{21}) = R(\beta ve_{22} + R(e_{21})),\\
0 = R(e_{21})R(1) = R(R(e_{21}) - \beta ve_{22}).
\end{gather*}
Subtracting these equalities, we get $2\beta R(ve_{22}) = 0$, whence $\beta\beta_2=0$ and $\beta = 0$. So, $R(1) = 0$. Since $\ker (R)$ is a~$\mathrm{Im}(R)$-bimodule, then it means that $ve_{22},e_{12}\in\ker (R)$. It is a~contradiction.

\medskip 2.d) $R(e_{12}) =  ve_{22} + \beta_1 ve_{12} - e_{12} = R(ve_{22}).$ It means that $\alpha + \gamma = 0$ by above. If $\alpha = 0$, then, as above, $\beta R(ve_{22}) = 0$, whence $\beta = 0$ and $R(1) = 0$, a~contradiction. Let us assume that $\alpha\neq 0$. Further, as we did before,
\begin{gather*}
L(ve_{12})\ni R(1)R(ve_{11}) = R(R(1)ve_{11} + R(ve_{11})) = R(\alpha e_{11} - \beta e_{12} + \alpha_3 ve_{22} + \gamma_3 e_{12}),\\
L(ve_{12})\ni R(ve_{11})R(1) = R(\alpha_3 ve_{22} + \gamma_3 e_{12} + \alpha e_{22} + \beta e_{12}).
\end{gather*}
Adding these equalities, we get $\alpha R(1) + 2\alpha_3 R(ve_{22}) + 2\gamma_3 R(e_{12})\in L(ve_{12})$, whence we have $\alpha^2 + 2\alpha_3 + 2\gamma_3 = 0$. Subtracting these equalities, we get $\alpha R(2e_{11}-1) - 2 \beta R(e_{12}) \in L(ve_{12})$. Therefore, $\alpha(2\alpha_5-\alpha) - 2\beta = 0$ and $\alpha(2\gamma_5 - \gamma) + 2\beta = 0$, whence $\alpha(\alpha_5+\gamma_5)=0$ and $\alpha_5+\gamma_5 = 0$.
Further,
\begin{gather*}
L(ve_{12}) \ni R(1)R(ve_{21}) = R(\alpha e_{21} - \beta e_{22} - \gamma ve_{11} + \alpha_4 ve_{22} + \gamma_4 e_{12}),\\
L(ve_{12}) \ni R(ve_{21})R(1) = R(\alpha_4 ve_{22} + \gamma_4 e_{12} - \alpha e_{21} - \beta e_{11} + \gamma ve_{11}).
\end{gather*}
Adding these equalities, we get $-\beta R(1) + 2\alpha_4 ve_{22} + 2\gamma_4 e_{12} \in L(ve_{12})$, whence $-\beta\alpha + 2\alpha_4+2\gamma_4 = 0$. Subtracting these equalities, we get $2\alpha e_{21} + \beta (2e_{11}-1) - 2\gamma ve_{11} \in L(ve_{12})$, whence $2\alpha\alpha_6 + \beta(2\alpha_5 - \alpha) - 2\gamma\alpha_3 = 0$, $2\alpha\gamma_6 + \beta(2\gamma_5 - \gamma) - 2\gamma\gamma_3 = 0$. Adding last two equalities, we get $\alpha(\alpha_6+\gamma_6 + \alpha_3+\gamma_3) = 0$, whence $\alpha_3+\gamma_3+\alpha_6+\gamma_6 = 0$. Further,
\begin{gather*}
L(ve_{12}) \ni R(1)R(e_{21}) = R(\beta ve_{22} + \gamma e_{11} + \alpha_6 ve_{22} + \gamma_6 e_{12}),\\
L(ve_{12}) \ni R(e_{21})R(1) = R(\alpha_6 ve_{22} + \gamma_6 e_{12} -\beta ve_{22} + \gamma e_{22}).
\end{gather*}
Adding these equalities, we get $\gamma R(1) + 2\alpha_6 ve_{22} + 2\gamma_6 e_{12} \in L(ve_{12})$, whence $\gamma\alpha + 2\alpha_6 + 2\gamma_6 = 0$. Subtracting these equalities, we get $\gamma(2e_{11}-1) + 2\beta ve_{22} \in L(ve_{12})$, whence $\gamma(2\alpha_5 - \alpha) + 2\beta = 0$. We have
\[L(ve_{12})\ni R(e_{21})R(e_{12}) = R(- \beta_1 ve_{22} + e_{22}),\]
whence $-\beta_1+\alpha-\alpha_5 = 0$. Further,
\[L(ve_{12})\ni R(e_{12})R(e_{21}) = R(\beta_1 ve_{22} - e_{11}),\]
whence $\beta_1 - \alpha_5 = 0$, $-\beta_1-\gamma_5 = 0$, i.e. $\alpha_5=-\gamma_5=\beta_1$. It means that $\alpha=2\alpha_5$ and $\beta = 0$. Then $\alpha_4+\gamma_4 = 0$ and $\alpha_6+\alpha_3 = 0$. Further,
\[L(ve_{12}) \ni R(ve_{21})R(e_{12}) = R(e_{21} - \beta_1e_{11} + ve_{11}),\]
whence $\alpha_6 - \beta_1\alpha_5 + \alpha_3 = 0$. It means that $\beta_1\alpha_5 = 0$, so $\alpha_5=\gamma_5=\beta_1=\alpha=\gamma=0$ and $R(1)=0$. We have a~contradiction.

\medskip 3) Let us assume that $\dim V = 1$, i.e. $V = L(ve_{12})$. Then $\dim\mathrm{Im}(R^2)=2$. Let
\begin{gather*}
R(1) = \alpha ve_{22} + \beta ve_{12} + \gamma e_{12},\\
R(e_{12})=\alpha_1 ve_{22} + \beta_1 ve_{12} + \gamma_1 e_{12},\\
R(ve_{22}) = \alpha_2 ve_{22} + \beta_2 ve_{12} + \gamma_2 e_{12}.
\end{gather*}
Then
\[0=R(1)R(1)=2R(R(1)) = 2R(\alpha ve_{22} + \gamma e_{12}) = 2\alpha R(ve_{22}) + 2\gamma R(e_{12}).\]
Since $\dim \mathrm{Im}(R^2) = 2$, then vectors $R(ve_{22})$ and $R(e_{12})$ are linearly independent. It means that $\alpha=\gamma=0$ and $R(1)=\beta ve_{12}$.

For any $x,y\in\mathbb{O}$ we have
\begin{gather*}
0 = R(1)R(x) = R(R(1)x + R(x)) = \beta R(ve_{12}x) + R^2(x),\\
0 = R(y)R(1) = R(R(y) + yR(1)) = R^2(y) + \beta R(yve_{12}).
\end{gather*}
For $x=e_{12},ve_{22}$ we obtain $R^2(e_{12})=R^2(ve_{22})=0$. Let $x=y=e_{21}$, then we have
\begin{gather*}
0 = \beta R(ve_{22}) + R^2(e_{21}),\\
0 = -\beta R(ve_{22}) + R^2(e_{21}).
\end{gather*}
Subtracting one of these equalities from another, we get $\beta R(ve_{22})=0$, whence $\beta = 0$ and $R(1)=0$. Then we have $R^2(x)=0$ from the equalities above for any $x\in\mathbb{O}$. We have a~contradiction. The lemma is proven.

\medskip\textbf{Corollary 3.} \textit{Let $R$ be a~Rota-Baxter operator on the split Cayley-Dickson algebra~$\mathbb{O}$ and $\mathrm{Im}(R)=Fe_{12} + Fve_{12} + Fve_{22}$. Then \[\ker (R) = L(e_{11},e_{12},e_{22},ve_{12},ve_{22}).\]
If a~field $F$ is quadratically closed then, up to conjugation by automorphism, antiautomorphism and up to multiplication by a~scalar, an operator $R$ acts in one of the following ways for some $\alpha,\beta,\gamma\in F$:
\begin{enumerate}
\item $R(e_{21})= e_{12},\quad R(ve_{11})=ve_{12},\quad R(ve_{21})=ve_{22} + ve_{12};$
\item $R(e_{21}) = e_{12},\quad R(ve_{11})=ve_{12},\quad R(ve_{21})=ve_{22};$
\item $R(e_{21})= e_{12},\quad R(ve_{11})=ve_{12},\quad R(ve_{21})=ve_{22} + e_{12},$
\item $R(e_{21})= \alpha e_{12},\quad R(ve_{11})=ve_{12},\quad R(ve_{21})=ve_{22} + ve_{12} + e_{12},$
where $\alpha\neq 0$.
\end{enumerate}
}
\textbf{Proof.} 1) A conjugation by the automorphism from the Proposition~8 with a~scalar $\frac{1}{\sqrt{\alpha}}$ and a~multiplication by a~scalar $\frac{1}{\alpha}$ give us the result.

2) A conjugation by the automorphism from the Proposition~8 with a~scalar $\frac{1}{\sqrt{\alpha}}$ and a~multiplication by a~scalar $\frac{1}{\alpha}$ give us the result. The corollary is proven.

\medskip \textbf{Lemma 6.} \textit{Let $R$ be the Rota-Baxter operator on the split Cayley-Dickson algebra $\mathbb{O}$ and $\mathrm{Im}(R)=Fe_{11} + Fve_{12} + Fve_{22}$. Then, up to the action of automorphism, antiautomorphism and multiplication by a~scalar, we can assume that $\ker (R) = L(e_{11},e_{12},e_{22},ve_{12}, ve_{22})$ and $R$ acts on $\mathbb{O}$ in one of the following ways:
\begin{enumerate}
    \item $R(e_{21}) =  ve_{12}, \quad R(ve_{11}) = e_{11},\quad R(ve_{21}) = ve_{22}$;
\item $R(e_{21}) =  \alpha e_{11},\quad R(ve_{11}) = \beta ve_{12},\quad R(ve_{21}) = ve_{22}$,
where $\alpha,\beta\neq 0$;
\item $R(e_{21}) =  -ve_{22},\quad R(ve_{11}) =  ve_{12},\quad R(ve_{21}) = e_{11}$.
\end{enumerate}
}

\textbf{Proof.} Let 
\begin{gather*}
    R(e_{11}) = \alpha_{1}ve_{22} + \beta_1 ve_{12} + \gamma_1 e_{11},\\
    R(ve_{12}) = \alpha_2 ve_{22} + \beta_2 ve_{12} + \gamma_2 e_{11},\\
    R(ve_{22}) = \alpha_3 ve_{22} + \beta_3 ve_{12} + \gamma_3 e_{11}.
\end{gather*}
Then
\[\gamma_1^2 e_{11} + \beta_1\gamma_1 ve_{12} + \gamma_1\alpha_1 ve_{22} = R(e_{11})R(e_{11}) = R(\beta_1 ve_{12} + 2\gamma_1 e_{11} + \alpha_1 ve_{22}),\]
whence 
    \[\alpha_1\alpha_3 + \beta_1\alpha_2 +\gamma_1\alpha_1 = \alpha_1\beta_3 + \beta_1\beta_2 + \gamma_1\beta_1 = \alpha_1\gamma_3 + \beta_1\gamma_2 + \gamma_1^2 = 0.\]
Further,
\[\gamma_2\gamma_3 e_{11} + \beta_2\gamma_3 ve_{12} + \gamma_2\alpha_3 ve_{22} = R(ve_{12})R(ve_{22}) = R(\gamma_2 ve_{22} + \gamma_3 ve_{12}),\]
where 
\[\gamma_3\alpha_2 = 0,\quad \beta_3\gamma_2 = 0,\quad \gamma_2\gamma_3 = 0.\]
An equality $R(e_{11})R(ve_{12}) = R(R(e_{11})ve_{12}+e_{11}R(ve_{12}))$ gives us
\[\alpha_2\beta_3 = 0,\quad \gamma_1\alpha_2 = \gamma_2\alpha_1 + \alpha_2\alpha_3.\]
An equality $R(ve_{22})R(e_{11}) = R(R(ve_{22})e_{11} + ve_{22}R(e_{11}))$ implies
\[\beta_3\gamma_1=\gamma_3\beta_1+\beta_3\beta_2.\]
An equality $R(e_{11})R(ve_{22}) = R(R(e_{11})ve_{22}+e_{11}R(ve_{22}))$ implies
\[\gamma_3(\gamma_1+\alpha_3) = 0,\quad \beta_3 (\gamma_1+\alpha_3) = 0,\quad \gamma_3\alpha_1 + \alpha_3^2 = 0. \]
Finally, an equality $R(ve_{12})R(e_{11}) = R(R(ve_{12})e_{11}+ve_{12}R(e_{11}))$ gives us
\[\gamma_2(\beta_2 +\gamma_1) = 0,\quad \gamma_2\beta_1 + \beta_2^2 = 0,\quad \alpha_2(\beta_2+\gamma_1) = 0.\]

1) Let $\gamma_2 \neq 0$. Then $\beta_3=\gamma_3 = 0$, whence $\alpha_3 = 0$ and $R(ve_{22}) = 0$.

2) Let $\beta_3\neq 0$. Then $\alpha_2 = \gamma_2 = 0$, whence $\beta_2 = 0$ and $R(ve_{12}) = 0$.

3) Let $\gamma_2 = \beta_3 = 0$. Then $\beta_2 = 0$. Let us consider two subcases.

3.a) Let $\gamma_1\neq 0$. Then $\beta_1 = \alpha_2 = 0$, whence $R(ve_{12}) = 0$.

3.b) Let $\gamma_1 = 0$. Then $\alpha_2\alpha_3 = 0$. If $\alpha_2 = 0$, then $R(ve_{12}) = 0$. If $\alpha_2 \neq 0$, then $\alpha_3 = \gamma_3 = 0$, whence $R(ve_{22}) = 0$.

So, in any case we have either $R(ve_{12}) = 0$, or $R(ve_{22}) = 0$. An antiautomorphism from Proposition~12 swaps the subspaces $L(ve_{12})$ and $L(ve_{22})$, but an element $e_{11}$ is a~fixed point. It allows us to assume that $R(ve_{22}) = 0$, i.e. $\alpha_3=\beta_3=\gamma_3=0$. Then the equalities above give us the following restrictions:
\begin{gather}
    \label{5.7}\beta_1\alpha_2 + \gamma_1 \alpha_1 = 0,\\
    \label{5.8}\gamma_1\alpha_2 = \gamma_2\alpha_1,\\
    \label{5.9}\beta_1\gamma_2 + \gamma_1^2 = 0,\\
    \label{5.10}\beta_1\gamma_2 + \beta_2^2 = 0,\\
    \label{5.11}\gamma_2(\beta_2+\gamma_1) = 0,\\
    \label{5.12}\beta_1(\beta_2+\gamma_1) = 0,\\
    \label{5.13}\alpha_2(\beta_2+\gamma_1) = 0.
\end{gather}
Let $\beta_2+\gamma_1 \neq 0$. Then by \eqref{5.11}--\eqref{5.13} we have $\alpha_2=\gamma_2=\beta_1 = 0$, whence $\beta_2 = \gamma_1 = 0$ by \eqref{5.9}--\eqref{5.10}, it is a~contradiction. Thus, $\beta_2 + \gamma_1 = 0$.

If $\gamma_1\neq 0$ then $\beta_1\neq 0$ by \eqref{5.9}--\eqref{5.10} and we have by \eqref{5.7}--\eqref{5.9}:
\begin{gather*}
    R(e_{11}) = \alpha_{1}ve_{22} + \beta_1 ve_{12} + \gamma_1 e_{11},\\
    R(ve_{12}) = -\frac{\gamma_1}{\beta_1}( \alpha_1ve_{22} + \beta_1 ve_{12} + \gamma_1 e_{11}).
\end{gather*}
It means that $e_{11}+\frac{\beta_1}{\gamma_1} ve_{12} \in\ker (R)$. But $\ker(R)$ is a~$Im(R)$-bimodule, so $e_{11}=e_{11}\cdot (e_{11}+\frac{\beta_1}{\gamma_1}ve_{12}) \in\ker (R)$, whence $\gamma_1 = 0$. It is a~contradiction, so $\gamma_1 = 0$. 

If $\gamma_2 \neq 0$, then $R(e_{11})=0$. But $Ker(R)$ is $Im(R)$-bimodule, so $ve_{12}=ve_{12}e_{11}\in Ker(R)$, it is a~contradiction. Thus, $\gamma_2=0$ and, up to multiplication by a~scalar, we have $\alpha_2\beta_1 = 0$ by \eqref{5.7}. 

If $\beta_1 = 0$ then we have $R(e_{11}) = \alpha_1 ve_{22}$, $R(ve_{12}) = \alpha_2 ve_{22}$. In this case if $\alpha_2\neq 0$, then $e_{11}-\frac{\alpha_1}{\alpha_2}ve_{12}\in\ker (R)$. But $Ker(R)$ is a~$Im(R)$-bimodule, so $e_{11}=e_{11}\cdot (e_{11}-\frac{\alpha_1}{\alpha_2}ve_{12}) \in\ker (R)$, whence $ve_{12} = ve_{12}\cdot e_{11}\in\ker (R)$ and $R^2 = 0$, a~contradiction. It means that in this case $\alpha_2 = 0$.

So, we have three possible cases
\begin{enumerate}
    \item $R(ve_{12})=R(ve_{22})=0$, $R(e_{11})=\alpha_1 ve_{22} + ve_{12}$;
    \item $R(ve_{12})=R(ve_{22}) = 0$, $R(e_{11}) = ve_{22}$;
    \item $R(ve_{12})=R(ve_{22})=R(e_{11})=0$.
\end{enumerate}
Let $x=\sum\alpha_{ij}e_{ij} + \sum\beta_{ij}e_{ij} \in \ker (R)$. Since $\ker (R)$ is a~$\mathrm{Im}(R)$-bimodule and $ve_{12},ve_{22}\in\ker (R)$, then $\ker (R)$ contains the following elements:
\begin{gather*} e_{11}xe_{11} = \alpha_{11}e_{11},\\
e_{11}x - \beta_{22}ve_{22} - e_{11}xe_{11} = \alpha_{12}e_{12} + \beta_{21}ve_{21},\\
ve_{22}(e_{11}x - \beta_{22}ve_{22} - e_{11}xe_{11}) -\alpha_{12}ve_{12} = \beta_{21}e_{21},\\
ve_{12}(e_{11}x - \beta_{22}ve_{22} - e_{11}xe_{11}) = -\beta_{21}e_{22},\\
(e_{11}x - \beta_{22}ve_{22} - e_{11}xe_{11})ve_{12} = -\beta_{21}e_{11},\\
x - e_{11}x - \beta_{12}ve_{12} = \alpha_{21}e_{21} + \alpha_{22}e_{22} + \beta_{11}ve_{11},\\
(x-e_{11}x)ve_{22} = \beta_{11}e_{22},\\
(x-e_{11}x)e_{11} - \beta_{12}ve_{12} = \alpha_{21}e_{21}+\beta_{11}ve_{11},\\
(x-e_{11}x)-(x-e_{11}x)e_{11} = \alpha_{22}e_{22},\\
ve_{22}((x-e_{11}x)e_{11}) = \beta_{11}e_{11},\\
ve_{12}((x-e_{11}x)e_{11}) -\alpha_{21}ve_{22} = -\beta_{11}e_{12}.
\end{gather*}
Let
\begin{gather*} 
R(ve_{11})=\alpha_4 ve_{22} + \beta_4 ve_{12}+\gamma_4 e_{11},\\
R(ve_{21})=\alpha_5 ve_{22} + \beta_5 ve_{12}+ \gamma_5 e_{11}.
\end{gather*}

Let us consider the three cases above.

\medskip (1) Let $R(ve_{12})=R(ve_{22})=0$, $R(e_{11})=\alpha_1 ve_{22} + ve_{12}$. Then $e_{11}\notin \ker (R)$, so $\alpha_{11}=\beta_{21}=\beta_{11} = 0$. Hence $\ker (R)=L(e_{12},e_{21},e_{22},ve_{12},ve_{22})$. Then
\[\gamma_4ve_{12} = R(e_{11})R(ve_{11}) = (\alpha_1+\gamma_4)\alpha_1 ve_{22} + (\alpha_1+\gamma_4) ve_{12},\]
whence $\alpha_1=0$ and $R(e_{11})=ve_{12}$. Further, 
\[0=R(ve_{21})R(e_{11})=R(\gamma_5 e_{11} + \beta_5 ve_{12} - e_{11}) = (\gamma_5-1)ve_{12},\]
whence $\gamma_5 = 1$. Further,
\begin{multline*}
\beta_4\gamma_5 ve_{12} + \gamma_4\gamma_5 e_{11} + \gamma_4\alpha_5 ve_{22}=R(ve_{11})R(ve_{21}) = (\gamma_4\alpha_5 + \gamma_5\alpha_4)ve_{22} + (\gamma_4\beta_5 + \gamma_5\beta_4)ve_{12} + 2\gamma_4\gamma_5 e_{11},
\end{multline*}
whence
\begin{gather*}
    \gamma_4\gamma_5 = 0,\quad \gamma_5\alpha_4 = 0,\quad \gamma_4\beta_5 = 0.
\end{gather*}
Since $\gamma_5=1$, then $\gamma_4=\alpha_4 = 0$ and $R(ve_{11})=\beta_4 ve_{12}$. But then $\dim(\mathrm{Im}(R)) < 3$, we have a~contradiction.

\medskip (2) Let $R(ve_{22}) = R(ve_{12}) = 0$, $R(e_{11}) = ve_{22}$. We have $e_{11}\notin\ker (R)$ and $\alpha_{11}=\beta_{21}=\beta_{11}=0$. It means that $\ker (R)=L(e_{12},e_{21},e_{22},ve_{12},ve_{22})$. Then
\[0=R(e_{11})R(ve_{11}) = R(e_{11} + \alpha_4 ve_{22} + \gamma_4 e_{11}) = (1+\gamma_4) ve_{22},\]
whence $\gamma_4 = -1$. Further,
\[0=R(e_{11})R(ve_{21}) = R(e_{21} + \gamma_5 e_{11} + \alpha_5 ve_{22}) = \gamma_5 ve_{22},\]
Whence $\gamma_5 = 0$. Finally,
\[-\alpha_5 ve_{22}=R(ve_{11})R(ve_{21}) = R(-ve_{21}) = -\alpha_5 ve_{22} -\beta_5 ve_{12},\]
whence $\beta_5 = 0$ and $R(ve_{21})=\alpha_5 ve_{22}$. But then we have $\dim(\mathrm{Im}(R)) < 3$, it is a~contradiction.

\medskip (3) Let $R(ve_{22})=R(ve_{12})=R(e_{11})=0$, i.e. $R^2 = 0$.
If $\beta_{21}\neq 0$, then $e_{12},e_{22},e_{21}\in\ker (R)$, it is a~contradiction with a~dimension of $\ker (R)$. Thus, $\beta_{21}=0$ and $\alpha_{12}e_{12}\in\ker (R)$. If $\beta_{11}\neq 0$, then $e_{12},e_{22}\in\ker (R)$ and $\alpha_{21}e_{21}+\beta_{11}ve_{11}\in\ker (R)$, i.e.
$\ker (R)=L(e_{11},e_{12},e_{22},ve_{12},ve_{22},\alpha_{21}e_{21}+\beta_{11}ve_{11})$, it is a~contradiction with a~dimension of $\ker (R)$. Thus, $\beta_{11}=0$ and $\alpha_{21}e_{21}\in\ker (R)$. It means that the set $\{e_{12},e_{21},e_{22}\}\cap\ker (R)$ contains two elements. 

Let
\begin{gather*}
    R(e_{12})=\alpha_6 ve_{22} + \beta_6 ve_{12}+\gamma_6 e_{11},\\
    R(e_{21})=\alpha_7 ve_{22} + \beta_7 ve_{12}+\gamma_7 e_{11},\\
    R(e_{22})=\alpha_8 ve_{22} + \beta_8 ve_{12}+\gamma_8 e_{11}.
\end{gather*}

(3.a) Let $R(e_{12})=R(e_{21})=0$. Then $R(e_{22})\neq 0$ and
\[\gamma_4\gamma_4 e_{11} + \gamma_4\alpha_4 ve_{22} + \gamma_4\beta_4 ve_{12}=R(ve_{11})R(ve_{11}) = \gamma_4^2 e_{11} + \gamma_4\alpha_4 ve_{22} + \gamma_4\beta_4 ve_{12} + \alpha_4R(e_{22}),
\]
whence $\alpha_4 R(e_{22})=0$, i.e. $\alpha_4=0$. Further,
\[\gamma_5\gamma_5 e_{11} + \gamma_5\alpha_5 ve_{22} + \gamma_5\beta_5 ve_{12}=R(ve_{21})R(ve_{21}) = \gamma_5^2 e_{11} + \gamma_5\alpha_5 ve_{22} + \gamma_5\beta_5 ve_{12} - \beta_5 R(e_{22}),\]
whence $\beta_5 R(e_{22})=0$, i.e. $\beta_5=0$. 

Since $R(1)R(1)=2R(R(1))=0$ and $R(e_{11})=0$, then $R(e_{22})R(e_{22})=2R(R(e_{22}))=0$. Then we have $\gamma_8=0$. Since $R(e_{11})=0$, then
\[R(e_{22})R(x)=R(1)R(x)=R(R(e_{22})x+R(x)) = R(R(e_{22})x).\]
For $x=ve_{11}$ we have
\[\gamma_4\beta_8 ve_{12} = R(-\beta_8 e_{12} + \alpha_8 e_{11}) = 0,\]
whence $\gamma_4\beta_8 = 0$. For $x=ve_{21}$ we have
\[\gamma_5\beta_8 ve_{12} = R(-\beta_8 e_{22} + \alpha_8 e_{21}) = -\beta_8^2 ve_{12} -\beta_8\alpha_8 ve_{22},\]
whence $\beta_8\alpha_8 = 0$ and $\beta_8(\gamma_5+\beta_8)=0$. If $\beta_8\neq 0$, then $\gamma_4=\alpha_8 = 0$. It means that $\dim(\mathrm{Im}(R)) < 3$, it is a~contradiction. So, $\beta_8 = 0$ and $\alpha_8\neq 0$. Similarly,
\[R(y)R(e_{22})=R(yR(e_{22})).\]
For $y=ve_{11}$ we get
\[\gamma_4\alpha_8 ve_{22} = R(\alpha_8 e_{22} + \beta_8 e_{12}) = \alpha_8^2 ve_{22} + \alpha_8\beta_8 ve_{12},\]
whence $\alpha_8(\gamma_4+\alpha_8) = 0$. But $\alpha_8\neq 0$, hence $\gamma_4+\alpha_8 = 0$. For $y=ve_{21}$ we get
\[\gamma_5\alpha_8 ve_{22} = R(-\alpha_8 e_{21} - \beta_8 e_{11}) = 0,\]
whence $\gamma_5\alpha_8 = 0$. But $\alpha_8\neq 0$, hence $\gamma_5 = 0$. But it means that $R(ve_{21})=\alpha_5 ve_{22}$ and $R(e_{22}) = \alpha_8 ve_{22}$. Hence $\dim(\mathrm{Im}(R)) < 3$, it is a~contradiction.

(3.b) Let $R(e_{12})=R(e_{22})=0$. Then
\[\gamma_7\gamma_4 e_{11} + \gamma_7\alpha_4 ve_{22} + \beta_7\gamma_4 ve_{12} = R(e_{21})R(ve_{11}) = \gamma_4\gamma_7 e_{11} + \gamma_4 \beta_7 ve_{12} + \gamma_4\alpha_7 ve_{22},\]
whence $\gamma_7\alpha_4 = \gamma_4\alpha_7$. Further,
\[R(e_{21})R(ve_{21}) = (\alpha_7+\gamma_5)R(e_{21}) + \gamma_7 R(ve_{21}),\]
whence 
\[(\alpha_7 + \gamma_5)\gamma_7  = 0,\quad (\alpha_7+\gamma_5)\alpha_7 = 0,\quad \alpha_7\beta_7 + \gamma_7\beta_5 = 0.\]
Further,
\[R(ve_{11})R(ve_{21}) = \alpha_4 R(e_{21}) + \gamma_4 R(ve_{21}) + \gamma_5 R(ve_{11}),\]
whence
\[ \gamma_4\gamma_5 + \alpha_4\gamma_7 = 0,\quad \alpha_4(\alpha_7 + \gamma_5) = 0,\quad \alpha_4\beta_7 + \gamma_4\beta_5 = 0. \]

If $\alpha_7+\gamma_5\neq 0$, then $\gamma_7=\alpha_7=\alpha_4 = 0$, whence $\gamma_4\gamma_5=0$. Since $\gamma_5\neq 0$, then $\gamma_4 = 0$ and $\dim(\mathrm{Im}(R)) < 3$ (because $R(e_{21}),R(ve_{11})\in L(ve_{12})$), it is a~contradiction. Therefore, $\alpha_7=-\gamma_5$. So, equalities above give us
\begin{gather}
    \label{5.14}\gamma_7\alpha_4 = -\gamma_4\gamma_5,\\
    \label{5.15}\gamma_5\beta_7 = \gamma_7\beta_5,\\
    \label{5.16}\alpha_4\beta_7 + \gamma_4\beta_5 = 0,\\
    \label{5.17}\alpha_7=-\gamma_5.
\end{gather}

(3.b.a) Let $\gamma_5 = 0$. Then we have $\alpha_7=0$, $\gamma_7\beta_5 = \gamma_7\alpha_4 = 0$ by \eqref{5.14}--\eqref{5.15},\eqref{5.17}.

(3.b.a.a) Let $\gamma_7 = 0$. Then $\beta_7\neq 0$, $\gamma_4\neq 0$. Up to multiplication by a~scalar we can assume that
\begin{gather*}
    R(e_{21}) =  ve_{12},\\
    R(ve_{11}) = \alpha_4 ve_{22} + \beta_4 ve_{12} + \gamma_4 e_{11},\\
    R(ve_{21}) = \alpha_5 ve_{22} + \beta_5 ve_{12},
\end{gather*}
where $\alpha_4 + \gamma_4\beta_5 = 0$ by \eqref{5.16}. A conjugation by an automorphism from Proposition~9 with a~scalar $\alpha = \beta_5$ allows us to assume that
\[R(e_{21}) =  ve_{12},\quad R(ve_{11}) = \beta_4 ve_{12} + \gamma_4 e_{11},\quad R(ve_{21}) = \alpha_5 ve_{22}.\]
The Proposition~8 with a~scalar $\alpha = \alpha_5$ and a~multiplication by a~scalar $\alpha_5$ allow us to assume that
\[R(e_{21}) =  ve_{12},\quad R(ve_{11}) = \beta_4 ve_{12} + \gamma_4 e_{11},\quad R(ve_{21}) = ve_{22}.\]
The Proposition~6 with a~scalar $\alpha = \gamma_4$ and a~multiplication by a~scalar $\frac{1}{\gamma_4}$ allow us to assume that
\[R(e_{21}) =  ve_{12},\quad R(ve_{11}) = \beta_4 ve_{12} + e_{11},\quad  R(ve_{21}) = ve_{22}.\]
The Proposition~13 with a~scalar $\beta_4$ allows us to assume that
\[R(e_{21}) =  ve_{12},\quad R(ve_{11}) = e_{11},\quad R(ve_{21}) = ve_{22}.\]

(3.b.a.b) Let $\gamma_7\neq 0$. Then we have $\beta_5 = \alpha_4 = 0$. Up to multiplication by a~scalar we can assume that
\[ R(e_{21}) =   \beta_7 ve_{12} + \gamma_7 e_{11},\quad R(ve_{11}) = \beta_4 ve_{12} + \gamma_4 e_{11},\quad  R(ve_{21}) = ve_{22}.\]

(3.b.a.b.a) If $\gamma_4 = 0$, then
\[R(e_{21}) =   \beta_7 ve_{12} + \gamma_7 e_{11},\quad R(ve_{11}) = \beta_4 ve_{12},\quad R(ve_{21}) = ve_{22}.\]
Proposition~15 with a~scalar $\frac{\beta_7}{\beta_4}$ allows us to assume that
\[R(e_{21}) =   \gamma_7 e_{11},\quad R(ve_{11}) = \beta_4 ve_{12},\quad R(ve_{21}) = ve_{22}.\]

(3.b.a.b.b) If $\gamma_4\neq 0$, then a~conjugation by an automorphism from the Proposition~8 with a~scalar $\frac{1}{\gamma_4}$ and a~multiplication by a~scalar $\frac{1}{\gamma_4^2}$ allow us to assume that
\[R(e_{21}) =   \beta_7 ve_{12} + \gamma_7 e_{11},\quad R(ve_{11}) = \beta_4 ve_{12} + e_{11},\quad R(ve_{21}) = ve_{22}.\]
The Proposition~13 with a~scalar $\frac{1}{\gamma_7}$ allows us to assume that
\[R(e_{21}) =   \beta_7 ve_{12} + \gamma_7 e_{11},\quad R(ve_{11}) = \beta_4' ve_{12},\quad R(ve_{21}) = ve_{22}.\]
We have a~case $\gamma_4 = 0$, it was considered in (3.b.a.b).

(3.b.b) Let $\gamma_5\neq 0$. By \eqref{5.14}--\eqref{5.17} we can assume that
\begin{gather*}
    R(e_{21}) =  -ve_{22} + \gamma_7\beta_5 ve_{12} + \gamma_7 e_{11},\\
    R(ve_{11}) = \alpha_4 ve_{22} + \beta_4 ve_{12} -\gamma_7\alpha_4 e_{11},\\
    R(ve_{21}) = \alpha_5 ve_{22} + \beta_5 ve_{12} + e_{11}.
\end{gather*}
Proposition~13 with a~scalar $-\alpha_4$ allows us to assume that 
\begin{gather*}
    R(e_{21}) =  -ve_{22} + \gamma_7\beta_5 ve_{12} + \gamma_7 e_{11},\\
    R(ve_{11}) =  \beta_4' ve_{12},\\
    R(ve_{21}) = \alpha_5 ve_{22} + \beta_5 ve_{12} + e_{11}.
\end{gather*}

(3.b.b.a) Let $\gamma_7 = 0$. Then we have
\[ R(e_{21}) =  -ve_{22},\quad R(ve_{11}) =  \beta_4' ve_{12},\quad R(ve_{21}) = \alpha_5 ve_{22} + \beta_5 ve_{12} + e_{11}.\]
Proposition~8 with a~scalar $\beta_4'$ and a~multiplication by a~scalar $\beta_4'$ allow us to assume that
\[R(e_{21}) =  -ve_{22},\quad R(ve_{11}) =  ve_{12},\quad R(ve_{21}) = \alpha_5 ve_{22} + \beta_5 ve_{12} + e_{11}.\]
Proposition~15 with a~scalar $\alpha = \beta_5$ allows us to assume that
\[R(e_{21}) =  -ve_{22},\quad R(ve_{11}) =  ve_{12},\quad R(ve_{21}) = \alpha_5 ve_{22} + e_{11}.\]
Proposition~9 with a~scalar $-\frac{\alpha_5}{2}$ allows us to assume that
\[R(e_{21}) =  -ve_{22},\quad R(ve_{11}) =  ve_{12},\quad R(ve_{21}) = e_{11}.\]

(3.b.b.b) Let $\gamma_7\neq 0$. Proposition~8 with a~scalar $\frac{1}{\gamma_7}$ and a~multiplication by a~scalar $\frac{1}{\gamma_7}$ allow us to assume that
\begin{gather*}
    R(e_{21}) =  -ve_{22} + \gamma_7\beta_5 ve_{12} + e_{11},\quad R(ve_{11}) =  \beta_4'' ve_{12},\\
    R(ve_{21}) = \alpha_5' ve_{22} + \beta_5' ve_{12} + e_{11}.
\end{gather*}
Proposition~9 with a~scalar $\alpha = 1$ allows us to assume that
\begin{gather*}
    R(e_{21}) =  \gamma_7\beta_5 ve_{12} + e_{11},\quad   R(ve_{11}) =  \beta_4'' ve_{12},\\
    R(ve_{21}) = (\alpha_5'+1) ve_{22} + (\beta_5'-\alpha\gamma_7\beta_5) ve_{12}.
\end{gather*}
So, we have a~case $\gamma_5 = 0$, it was considered in (3.b.a).

(3.c) Let $R(e_{21})=R(e_{22})=0$. A conjugation by an antiautomorphism from the Proposition~12 allows us to assume that $R(e_{12})=0$, i.e. we are in the previous case. The lemma is proven.

\medskip \textbf{Corollary 4.} \textit{Let $R$ be the Rota-Baxter operator on the split Cayley-Dickson algebra $\mathbb{O}$ and $\mathrm{Im}(R)=Fe_{11} + Fve_{12} + Fve_{22}$. If a~field $F$ is quadratically closed then, up to the action of automorphism, antiautomorphism and multiplication by a~scalar, we can assume that $\ker (R) = L(e_{11},e_{12},e_{22},ve_{12}, ve_{22})$ and $R$ acts on $\mathbb{O}$ in one of the following ways:
\begin{enumerate}
    \item $R(e_{21}) =  ve_{12}, \quad R(ve_{11}) = e_{11},\quad R(ve_{21}) = ve_{22}$;
\item $R(e_{21}) =  e_{11},\quad R(ve_{11}) = \alpha ve_{12},\quad R(ve_{21}) = ve_{22}$,
where $\alpha\neq 0$;
\item $R(e_{21}) =  -ve_{22},\quad R(ve_{11}) =  ve_{12},\quad R(ve_{21}) = e_{11}$.
\end{enumerate}
}

\textbf{Proof.} 2) Proposition~7 with a~scalar $\frac{1}{\sqrt{\gamma_7}}$ and a~multiplication by a~scalar $\frac{1}{\sqrt{\gamma_7}}$ allows us to assume that
\[R(e_{21}) =   e_{11},\quad R(ve_{11}) = \beta ve_{12},\quad R(ve_{21}) = ve_{22}.\]
The corollary is proven.
 
\section{RB-Operators with four-dimensional image}

In \cite{Octonions} it was proved that there is only on four-dimensional non-unital subalgebra $B$ in $\mathbb{O}$, up to action of automorphism and antiautomorphism: $Fe_{11}+Fve_{12}+Fve_{11}+Fve_{12}$. Let us describe the Rota-Baxter operators of zero weight on $\mathbb{O}$ with this image.

\medskip \textbf{Lemma 7.} \textit{Let $R$ be the Rota-Baxter operator on the split Cayley-Dickson algebra $\mathbb{O}$ and $\mathrm{Im}(R)=Fe_{11}+Fe_{12} + Fve_{11} + Fve_{12}$. Then, up to the action of automorphism, antiautomorphism and multiplication by a~scalar, $R$ acts on $\mathbb{O}$ in one of the following ways for some $\alpha\in F$ (an operator $R$ is zero on unspecified basic elements $e_{ij}$, $ve_{ij}$):
\begin{enumerate}
    \item  $R(e_{11}) =  e_{12},\quad R(e_{21}) = - e_{11} ,\quad R(ve_{21}) =  - ve_{11},\quad R(ve_{22}) =  - ve_{12}$;
    \item $R(e_{11}) =  e_{12},\, R(e_{21}) = - e_{11}  -  ve_{11},\, R(ve_{21}) =  - ve_{11},\, R(ve_{22}) =  e_{12} - ve_{12}$;
    \item $R(e_{11}) = ve_{12},\quad R(e_{21}) =  -\alpha ve_{11},\quad R(ve_{21}) = e_{11},\quad R(ve_{22}) = \alpha e_{12}$,
     $\alpha\neq 0$;
    \item $R(e_{11}) = -e_{12} + ve_{12},\, R(e_{21}) =  e_{11} - \alpha ve_{11},\, R(ve_{21}) = e_{11} + ve_{11}$, $R(ve_{22}) = \alpha e_{12} + ve_{12}$, $\alpha\neq -1$;
    \item $R(e_{22}) = ve_{12},\quad R(e_{21}) =  -\alpha ve_{11},\quad  R(ve_{21}) =  - e_{11},\, R(ve_{22}) = \alpha e_{12}$, $\alpha\neq 0$.
\end{enumerate}
}

\textbf{Proof.} Let $x=\sum\alpha_{ij}e_{ij} + \sum\beta_{ij}ve_{ij}\in\ker (R)$. Then $\ker (R)$ contains elements
\begin{gather*}e_{11}xe_{11} = \alpha_{11}e_{11},\\
e_{11}x - e_{11}xe_{11} = \alpha_{12}e_{12} + \beta_{21}ve_{21} + \beta_{22}ve_{22},\\
(e_{11}x - e_{11}xe_{11})ve_{12} = -\beta_{21}e_{11},\\
ve_{12}(e_{11}x - e_{11}xe_{11}) = -\beta_{21}e_{22},\\
(xe_{11} - e_{11}xe_{11})ve_{11} =  \beta_{22}e_{11},\\
ve_{11}(e_{11}x - e_{11}xe_{11}) = \beta_{22}e_{22},\\
xe_{11}-e_{11}xe_{11} = \alpha_{21}e_{21} + \beta_{11}ve_{11} + \beta_{12}ve_{12},\\
e_{12}(xe_{11}-e_{11}xe_{11}) = \alpha_{21}e_{11},\\
(xe_{11}-e_{11}xe_{11})e_{12} = \alpha_{21}e_{22}.
\end{gather*}

If $\beta_{21}\neq 0$, then $e_{11},e_{22}\in\ker (R)$, whence $ve_{11}\cdot e_{11} = ve_{11}\in\ker (R)$, $ve_{12}\cdot e_{11} = ve_{12}\in\ker (R)$ and $e_{11}e_{12} = e_{12} \in\ker (R)$. It means that $\dim\ker (R) > 4$, it is a~contradiction. Therefore, $\beta_{21} = 0$. Similarly, $\alpha_{21} = \beta_{22} = 0$.

If $\alpha_{22}\neq 0$, then $e_{22}\in\ker (R)$, whence $e_{12}e_{22}=e_{12}\in\ker (R)$, $e_{22}ve_{12} = ve_{12}\in\ker (R)$, $e_{22}ve_{11}=ve_{11}\in\ker (R)$ and $\ker (R) = L(e_{12},e_{22},ve_{11},ve_{12})$. If $\alpha_{22}=0$, then $\ker (R) = L(e_{11},e_{12},ve_{11},ve_{12})$. Let
\begin{gather*}
    R(1) = \alpha_{1}e_{11} + \beta_1 e_{12} + \gamma_1 ve_{11} + \delta_1 ve_{12},\\
    R(e_{21}) = \alpha_{2}e_{11} + \beta_2 e_{12} + \gamma_2 ve_{11} + \delta_2 ve_{12},\\
    R(ve_{21}) = \alpha_{3}e_{11} + \beta_3 e_{12} + \gamma_3 ve_{11} + \delta_3 ve_{12},\\
    R(ve_{22}) = \alpha_{4}e_{11} + \beta_4 e_{12} + \gamma_4 ve_{11} + \delta_4 ve_{12}.
\end{gather*}

\medskip So, $e_{12},ve_{11},ve_{12} \in\ker (R)$ and a~set $\{e_{11},e_{22}\}\cap\ker (R)$ contains precisely one element. Then 
\[\alpha_1R(1)=R(1)R(1) = 2R(R(1)) = 2\alpha_{1}R(e_{11}),\]
whence either $\alpha_1^2 = 2\alpha_1^2$ (if $e_{22}\in\ker (R)$), or $\alpha_1^2=0$ (if $e_{11}\in\ker (R)$). In any case $\alpha_1=0$.  Further,
\[\alpha_2R(e_{21}) = R(e_{21})R(e_{21}) = R(\beta_2\cdot 1 + \alpha_2 e_{21}),\]
whence $\beta_2 R(1)=0$, so $\beta_2=0$. Further,
\[\alpha_3 R(ve_{21}) = R(ve_{21})R(ve_{21}) = R(\alpha_3 ve_{21} - \delta_3\cdot 1),\]
whence $\delta_3 R(1)=0$, so $\delta_3 = 0$. Further,
\[\alpha_4 R(ve_{22}) = R(ve_{22})R(ve_{22}) = R(\alpha_4 ve_{22} + \gamma_4\cdot 1),\]
whence $\gamma_4 R(1) = 0$, so $\gamma_4 = 0$. 

Let us notice that $R(1),R(e_{21}),R(ve_{21}),R(ve_{22})\neq 0$, because otherwise a~dimension of image is less than four.

We have two possible cases. 

\medskip 1) $\ker (R) = L(e_{22},e_{12},ve_{11},ve_{12})$. Then
\[\alpha_2\alpha_3 e_{11} + (-\gamma_3\delta_2 + \alpha_2\beta_3)e_{12} + \gamma_2\alpha_3 ve_{11} + \delta_2\alpha_3 ve_{12} = R(e_{21})R(ve_{21}) = R(\alpha_3 e_{21} + (\alpha_2-\gamma_3) ve_{21}),\]
whence 
\begin{gather}
 \label{6.1}\alpha_3(\alpha_2 - \gamma_3)  = 0,\\
 \label{6.3}\gamma_3(\delta_2 - \beta_3) = 0.
\end{gather}
Further,
\[(\alpha_2\beta_1 + \gamma_2\delta_1 - \delta_2\gamma_1) e_{12} = R(e_{21})R(e_{11}) = R(\alpha_2 e_{11} - \gamma_1 ve_{21} - \delta_1 ve_{22}), \]
whence
\begin{gather}
    \label{6.4}-\gamma_1\alpha_3 = \delta_1\alpha_4,\\
    \label{6.5}\gamma_1 (\alpha_2 - \gamma_3) = 0,\\    
    \label{6.6}\delta_1(\alpha_2-\delta_4) = 0,\\
    \label{6.7}\gamma_2\delta_1 - \delta_2\gamma_1 = - \gamma_1\beta_3 - \delta_1\beta_4.
\end{gather}
Further,
\[(\alpha_3\beta_1 + \gamma_3\delta_1) e_{12} = R(ve_{21})R(e_{11}) = R(\alpha_3 e_{11} -\delta_1 e_{11}), \]
whence
\begin{gather}
    \label{6.8}(\alpha_3-\delta_1)\gamma_1 = 0,\\
    \label{6.9}(\alpha_3-\delta_1)\delta_1 = 0,\\   
    \label{6.10}\delta_1 (\gamma_3+\beta_1) = 0.
\end{gather}
Further,
\[(\alpha_4\beta_1 - \delta_4\gamma_1) e_{12} = R(ve_{22})R(e_{11}) = R(\alpha_4 e_{11} + \gamma_1 e_{11}),\]
whence 
\begin{gather}
    \label{6.11}(\alpha_4+\gamma_1)\gamma_1 = 0,\\   
    \label{6.13}\gamma_1(\beta_1+\delta_4) = 0.
\end{gather}
Further,
\[(\gamma_1\delta_2-\gamma_2\delta_1)e_{12} + \gamma_1\alpha_2 ve_{11} + \delta_1\alpha_2 ve_{12} = R(e_{11})R(e_{21}) = R(\beta_1 e_{11} + \gamma_1 ve_{21} + \delta_1 ve_{22} + \alpha_2 e_{11}),\]
whence 
\begin{gather}
    \label{6.14}\gamma_1\delta_2-\gamma_2\delta_1 = (\beta_1+\alpha_2)\beta_1 + \gamma_1\beta_3 + \delta_1\beta_4,\\
    \label{6.16}\delta_1(\beta_1 + \delta_4) = 0.
\end{gather}
From \eqref{6.7} and \eqref{6.14} we have
\begin{equation}
    (\beta_1+\alpha_2)\beta_1 = 0.\label{6.17}
\end{equation}
Further,
\[-\gamma_3\delta_1 e_{12} + \gamma_1\alpha_3 ve_{11} + \delta_1\alpha_3 ve_{12} = R(e_{11})R(ve_{21}) = R(\alpha_3 e_{11}),\]
whence
\begin{equation}
    -\gamma_3\delta_1 = \alpha_3\beta_1.\label{6.18}
\end{equation}
Further,
\[\gamma_1\delta_4 e_{12} + \gamma_1\alpha_4 ve_{11} + \delta_1\alpha_4 ve_{12} = R(e_{11})R(ve_{22}) = R(\alpha_4 e_{11}),\]
whence 
\begin{gather}
    \gamma_1\delta_4 = \alpha_4\beta_1.\label{6.19}
\end{gather}
Further,
\[\alpha_3\alpha_4 e_{11} + (\alpha_3\beta_4 + \gamma_3\delta_4) e_{12} + \gamma_3\alpha_4 ve_{11} = R(ve_{21})R(ve_{22}) = R(\alpha_3 ve_{22} - \delta_4 e_{11}),\]
whence
\begin{equation}
\label{6.20}\delta_4(\gamma_3 + \beta_1) = 0.\\
\end{equation}
Further,
\[\alpha_3\alpha_4 e_{11} + (\alpha_4\beta_3-\delta_4\gamma_3) e_{12} + \delta_4\alpha_3 ve_{12} = R(ve_{22})R(ve_{21}) = R(\alpha_4 ve_{21} + \gamma_3 e_{11}),\]
whence
\begin{equation}
    \label{6.23}\gamma_3(\delta_4 + \beta_1) = 0.
\end{equation}
Further,
\[\alpha_2\alpha_4 e_{11} + (\alpha_2\beta_4+\gamma_2\delta_4) e_{12} + \gamma_2\alpha_4 ve_{11} + \delta_2\alpha_4 ve_{12} = R(e_{21})R(ve_{22}) = R(\alpha_2 ve_{22} + \alpha_4 e_{21} - \delta_4 ve_{22}),\]
whence 
\begin{gather}
    \label{6.27}\delta_4(\gamma_2 + \beta_4) = 0,\\
    \label{6.28}\delta_4(\alpha_2 - \delta_4) = 0.
\end{gather}
Further,
\[\alpha_4\alpha_2 e_{11} - \delta_4\gamma_2 e_{12} + \delta_4\alpha_2 ve_{12} = R(ve_{22})R(e_{21}) = R(\beta_4 e_{11} + \delta_4 ve_{22} + \gamma_2 e_{11}),\]
whence
\begin{equation}
    \label{6.30}(\beta_4+\gamma_2)\gamma_1 = 0.
\end{equation}
 Further,
\[\alpha_3\alpha_2 e_{11} + \gamma_3\delta_2 e_{12} + \gamma_3\alpha_2 ve_{11} = R(ve_{21})R(e_{21}) = R(\beta_3 e_{11} + \gamma_3 ve_{21} - \delta_2 e_{11}),\]
whence
\begin{equation}
    \label{6.33}\gamma_3\alpha_2 = (\beta_3-\delta_2)\delta_1 + \gamma_3^2.
\end{equation}

1.a) Let $\gamma_1=\delta_1 = 0$. Then $\beta_1\neq 0$ (because $R(1)\neq 0$), so from \eqref{6.17} we have $\beta_1=-\alpha_2$, from \eqref{6.18} and $\eqref{6.19}$ we have $\alpha_3 = \alpha_4 = 0$.

1.a.a) Let $\delta_4 = 0$. Then $\gamma_3 = 0$ from \eqref{6.23}, and $\dim(\mathrm{Im}(R)) < 4$, it is a~contradiction.

1.a.b) Let $\delta_4\neq 0$. Then we have $\gamma_3=-\beta_1$ from \eqref{6.20}, $\gamma_2=-\beta_4$ from \eqref{6.27}, $\alpha_2=\delta_4$ from \eqref{6.28}. Since $\gamma_3=-\beta_1\neq 0$ then $\delta_2=\beta_3$ from \eqref{6.3}. So, after multiplication by a~scalar $\frac{1}{\beta_1}$ we have
\begin{gather*}
    R(e_{11}) =  e_{12} ,\\
    R(e_{21}) = - e_{11}  - \beta_4' ve_{11} + \beta_3' ve_{12},\\
    R(ve_{21}) =  \beta_3' e_{12} - ve_{11},\\
    R(ve_{22}) =  \beta_4' e_{12} - ve_{12},
\end{gather*}
where $\beta_3' = \frac{\beta_3}{\beta_1}$, $\beta_4' = \frac{\beta_4}{\beta_1}$. Proposition~15 with a~scalar $\alpha = -\frac{\beta_3'}{2}$ allows us to assume that $\beta_3'=0$.

1.a.b.a) Let $\beta_4' = 0$. Then we have a~case (1) from the statement of the lemma.

1.a.b.b)  Let $\beta_4'\neq 0$. Proposition~7 with a~scalar $\frac{1}{\beta_4'}$ and a~multiplication by a~scalar $\beta_4'$ allow us to assume that we have a~case (2) from the statement of the lemma.

1.b) Let $\gamma_1 = 0$, $\delta_1\neq 0$. Then we have $\alpha_4 = 0$ from \eqref{6.4}, $\gamma_2=-\beta_4$ from \eqref{6.7}, $\alpha_3=\delta_1\neq 0$ from \eqref{6.9}, $\alpha_2=\gamma_3$ from \eqref{6.1}, $\gamma_3=-\beta_1$ from \eqref{6.10}, $\beta_1=-\delta_4$ from \eqref{6.16}, 

1.b.a) Let $\delta_4 = 0$. Then we have $\alpha_2 = \beta_1 = \gamma_3 = 0$ from \eqref{6.1}, \eqref{6.10} and \eqref{6.16}. So, $\beta_3=\delta_2$ from \eqref{6.33}. We have
\begin{gather*}
    R(e_{11}) = \alpha_3 ve_{12},\\
    R(e_{21}) =  -\beta_4 ve_{11} + \beta_3 ve_{12},\\
    R(ve_{21}) = \alpha_{3}e_{11} + \beta_3 e_{12} ,\\
    R(ve_{22}) = \beta_4 e_{12}.
\end{gather*}
We can see that $\beta_4\neq 0$. Then Proposition~2 with a~scalar $\alpha = \frac{\beta_3}{\beta_4}$ allows us to assume that
\[R(e_{11}) = \alpha_3 ve_{12},\quad R(e_{21}) =  -\beta_4 ve_{11},\quad R(ve_{21}) = \alpha_{3}e_{11},\quad R(ve_{22}) = \beta_4 e_{12}.\]
Proposition~6 with a~scalar $\alpha = \frac{1}{\alpha_3}$ allows us to assume that we have a~case (3) from the statement of the lemma.

1.b.b) Let $\delta_4\neq 0$. Since $\alpha_3=\delta_1\neq 0$, then we have $\alpha_2=\gamma_3$ from \eqref{6.1}. Since $\gamma_3\neq 0$, then we have $\delta_2=\beta_3$ from \eqref{6.3}. Proposition~15 with a~scalar $\alpha = 2\frac{\delta_2}{\delta_4}$ allows us to assume that $\delta_2 = \beta_3 = 0$. So, after multiplication by a~scalar $\frac{1}{\delta_1}$ we have
\begin{gather*}
    R(e_{11}) = -e_{12} + \delta_1' ve_{12},\\
    R(e_{21}) =  e_{11} - \beta_4' ve_{11},\\
    R(ve_{21}) = \delta_1' e_{11} + ve_{11},\\
    R(ve_{22}) = \beta_4' e_{12} + ve_{12},
\end{gather*}
where $\delta_1' = \frac{\delta_1}{\delta_4}$, $\beta_4' = \frac{\beta_4}{\delta_4}$. Proposition~7 with a~scalar $\delta_1'$ and a~multiplication by a~scalar $\frac{1}{\delta_1'}$ allow us to assume that we have a~case (4) from the statement of the lemma.

1.c) Let $\gamma_1 \neq 0$, $\delta_1 = 0$. The automorphism from Proposition~4 allows us to assume that this case is equivalent to the case (1.b)

1.d) Let $\gamma_1,\delta_1\neq 0$. Proposition~2 with a~scalar $\alpha = -\frac{\delta_1}{\gamma_1}$ allows us to assume that $\delta_1 = 0$, so we have the case 1.c).

\medskip 2) $\ker (R) = L(e_{11},e_{12},ve_{11},ve_{12})$, i.e. $R^2=0$. Then we have
\[\gamma_1\alpha_3 ve_{11} -\gamma_3\delta_1 e_{12} + \delta_1\alpha_3 ve_{12}=R(e_{22})R(ve_{21}) = -\delta_1R(e_{22}),\]
whence 
\begin{gather}
\label{6.024}\delta_1(\alpha_3 + \delta_1) = 0,\\
\label{6.025}\delta_1(\gamma_3-\beta_1) = 0.
\end{gather}
Further,
\[\gamma_1\alpha_4 ve_{11} + \gamma_1\delta_4 e_{12} + \delta_1\alpha_4 ve_{12}=R(e_{22})R(ve_{22}) = \gamma_1 R(e_{22}),\]
whence
\begin{equation}
\label{6.027}\delta_1(\alpha_4-\gamma_1) = 0.
\end{equation}
Further, 
\[\gamma_1\alpha_2 ve_{11} + (\gamma_1\delta_2 - \gamma_2\delta_1)e_{12} + \delta_1\alpha_2 ve_{12} = R(e_{22})R(e_{21}) = R(\gamma_1 ve_{21} + \delta_1 ve_{22}),\]
whence 
\begin{gather}
    \label{6.029}\gamma_1(\delta_2-\beta_3) = \delta_1(\beta_4+\gamma_2),\\
    \label{6.031}\delta_1(\alpha_2 - \delta_4) = 0.
\end{gather}
Further,
\[(\alpha_3\beta_1 + \gamma_3\delta_1)e_{12} = R(ve_{21})R(e_{22}) = 0,\]
whence
\begin{equation}
    \label{6.032}\alpha_3\beta_1 = -\gamma_3\delta_1.
\end{equation}
Further,
\[(\alpha_4\beta_1 - \delta_4\gamma_1) e_{12}   = R(ve_{22})R(e_{22}) = 0,\]
whence
\begin{equation}
    \label{6.033}\alpha_4\beta_1 = \delta_4\gamma_1.
\end{equation}
Further,
\[(\alpha_2\beta_1 + \gamma_2\delta_1 - \delta_2\gamma_1) e_{12} = R(e_{21})R(e_{22}) = R(\beta_1 e_{22} - \gamma_1 ve_{21} - \delta_1 ve_{22}),\]
whence
\begin{gather}
    \label{6.036}\delta_1(\beta_1 - \delta_4) = 0,\\
    \label{6.037}(\alpha_2-\beta_1)\beta_1 = 0.
\end{gather}
Further, 
\[
\alpha_2\alpha_4 e_{11} + (\alpha_2\beta_4 + \gamma_2\delta_4) e_{12} + \gamma_2\alpha_4 ve_{11} + \delta_2\alpha_4 ve_{12} = R(e_{21})R(ve_{22}) = R(\alpha_2 ve_{22} + \gamma_2 e_{22} + \alpha_4 e_{21} + \beta_4 e_{22} - \delta_4 ve_{22}),
\]
whence
\begin{equation}
    \label{6.041}\alpha_4\delta_4 + (\gamma_2+\beta_4)\delta_1 - \delta_4^2 = 0.
\end{equation}
Further,
\[\alpha_2\alpha_3 e_{11} + (\alpha_2\beta_3 - \delta_2\gamma_3) e_{12} + \gamma_2\alpha_3 ve_{11} + \delta_2\alpha_3 ve_{12}  = R(e_{21})R(ve_{21}) = R(\alpha_2 ve_{21} - \delta_2 e_{22} + \alpha_3 e_{21} + \beta_3 e_{22} - \gamma_3 ve_{21}),
\]
whence
\begin{equation}
    \label{6.057}(\beta_3-\delta_2)\delta_1 = 0.
\end{equation}

2.a) Let $\gamma_1 = \delta_1 = 0$. Then $\beta_1\neq 0$, so $\alpha_4=\alpha_3 = 0$ by \eqref{6.032} and \eqref{6.033}, $\alpha_2=\beta_1$ by \eqref{6.037}, $\delta_4 = 0$ from \eqref{6.041}. So, the dimension of the image is less than 4 ($R(e_{22},R(ve_{22})\in L(e_{12}$), it is a~contradiction. 

2.b) Let $\gamma_1 = 0, \delta_1\neq 0$. Then we have $\alpha_3 = -\delta_1$ from \eqref{6.024}, $\gamma_3 = \beta_1$ from \eqref{6.025}, $\alpha_4 = \gamma_1 = 0$ from \eqref{6.027}, $\beta_4 = -\gamma_2$ from \eqref{6.029}, $\alpha_2 = \delta_4$ from \eqref{6.031}, $\beta_1 = \delta_4$ from \eqref{6.036}, $\delta_2 = \beta_3$ from \eqref{6.057}, $\delta_4 = 0$ from \eqref{6.041}. After multiplication by a~scalar $\frac{1}{\delta_1}$, we can assume that
\begin{gather*}
    R(e_{22}) = ve_{12},\quad R(e_{21}) =  -\beta_4 ve_{11} + \delta_2 ve_{12},\\
    R(ve_{21}) =  - e_{11} + \delta_2 e_{12},\quad R(ve_{22}) = \beta_4 e_{12}.
\end{gather*}
Proposition~2 with a~scalar $\frac{\delta_2}{\beta_4}$ allows us to assume that we have a~case (5) from the statement of the lemma.

2.c) Let $\gamma_1 \neq 0$, $\delta_1 = 0$. The automorphism from Proposition~4 allows us to assume that this case is equivalent to the case 2.b).

2.d) Let $\gamma_1\neq 0$, $\delta_1\neq 0$. Proposition~2 with a~scalar $\alpha = -\frac{\delta_1}{\gamma_1}$ allows us to assume that $\delta_1 = 0$, so we have the case 2.c). The lemma is proven.

\medskip \textbf{Corollary 5.} \textit{Let $R$ be the Rota-Baxter operator on the split Cayley-Dickson algebra $\mathbb{O}$ and $\mathrm{Im}(R)=Fe_{11}+Fe_{12} + Fve_{11} + Fve_{12}$. If a~field $F$ is quadratically closed then, up to the action of automorphism, antiautomorphism and multiplication by a~scalar, $R$ acts on $\mathbb{O}$ in one of the following ways for some $\alpha\in F$ (an operator $R$ is zero on unspecified basic elements $e_{ij}$, $ve_{ij}$):
\begin{enumerate}
    \item  $R(e_{11}) =  e_{12},\quad R(e_{21}) = - e_{11},\quad R(ve_{21}) =  - ve_{11},\quad R(ve_{22}) =  - ve_{12}$
    \item $R(e_{11}) =  e_{12},\, R(e_{21}) = - e_{11}  + ve_{12},\, R(ve_{21}) =  e_{12} - ve_{11},\, R(ve_{22}) =  - ve_{12}$;
    \item $R(e_{11}) = ve_{12},\quad R(e_{21}) =  - ve_{11},\quad R(ve_{21}) = e_{11},\quad R(ve_{22}) =  e_{12}$;
    \item $R(e_{11}) = -e_{12} + ve_{12}, \quad R(e_{21}) =  e_{11} - \alpha ve_{11}, \quad R(ve_{21}) = e_{11} + ve_{11}$,\quad $R(ve_{22}) = \alpha e_{12} + ve_{12},$
where $\alpha\neq -1$;
    \item $R(e_{22}) = ve_{12},\quad R(e_{21}) =  - ve_{11},\quad R(ve_{21}) =  - e_{11},\quad R(ve_{22}) =  e_{12}.$
\end{enumerate}
}

\textbf{Proof.} Let us consider cases from Lemma~7.

3) and 5) Proposition~7 with a~scalar $\frac{1}{\sqrt{\alpha}}$ gives us the result.

The corollary is proven.

\section{Main Theorem}

We are ready to formulate the main result.

\medskip\textbf{Theorem 1.} \textit{Let $R$ be the Rota-Baxter operator on the split Cayley-Dickson algebra $\mathbb{O}$. Then, up to the action of automorphism, antiautomorphism and multiplication by a~scalar, $R$ acts on $\mathbb{O}$ in one of the following ways for some $\alpha,\beta\in F$ (an operator $R$ is zero on unspecified basic elements $e_{ij}$, $ve_{ij}$):
\begin{enumerate}
    \item $R(e_{21}) = e_{12}$;
    \item $R(ve_{22}) = e_{12}$;
    \item $R(e_{21}) = e_{11}$;
    \item $R(e_{21})=e_{11}$, $R(e_{22})=e_{12}$;
    \item $R(e_{21})=-e_{11}$, $R(e_{11})=e_{12}$;
    \item $R(e_{21})=e_{11}$, $R(ve_{21})= e_{12}$;
    \item $R(ve_{11})=\alpha e_{11}$, $R(ve_{21})= e_{12}$, $\alpha\neq 0$;
    \item $R(ve_{11}) = e_{12}, R(ve_{21}) = e_{11}$;
    \item $R(ve_{21}) = \alpha e_{11}, R(ve_{22}) = e_{12}$, $\alpha\neq 0$;
    \item $R(ve_{11}) = ve_{22}$, $R(ve_{21}) = ve_{22} + \alpha ve_{12}$, $\alpha\neq 0$;
    \item $R(e_{21}) = ve_{12}$, $R(ve_{21}) = ve_{22}$;
    \item $R(e_{11}) = R(e_{12}) = -R(e_{21}) = - R(e_{22}) = ve_{22} + ve_{12}, R(ve_{11}) = - ve_{12}, R(ve_{21}) = ve_{12}$;
    \item $R(ve_{11}) = ve_{12}$, $R(ve_{21}) = ve_{22} + ve_{12}$;
    \item $R(ve_{11}) =  ve_{12}$, $R(ve_{21}) = \alpha ve_{22}$, $\alpha\neq 0$;
    \item $R(ve_{11}) = ve_{22}$, $R(ve_{21}) = \alpha ve_{12}$, $\alpha\neq 0$;
    \item $R(e_{21})=\alpha e_{12}$, $R(ve_{11})=ve_{12}$, $R(ve_{21})=ve_{22} + ve_{12}$, $\alpha\neq 0$;
    \item $R(e_{21}) = \alpha e_{12}$, $R(ve_{11})=ve_{12}$, $R(ve_{21})=ve_{22}$, $\alpha\neq 0$;
    \item $R(e_{21})= e_{12}$, $R(ve_{11})=ve_{12}$, $R(ve_{21})=ve_{22} + e_{12}$;
    \item $R(e_{21})=\alpha e_{12}$, $R(ve_{11})=ve_{12}$, $R(ve_{21})=ve_{22} + ve_{12} + e_{12}$, $\alpha\neq 0$;
    \item $R(e_{21}) =  ve_{12}$, $R(ve_{11}) = e_{11}$, $R(ve_{21}) = ve_{22}$;
    \item $R(e_{21}) =  \alpha e_{11}$, $R(ve_{11}) = \beta ve_{12}$, $R(ve_{21}) = ve_{22}$, $\alpha,\beta\neq 0$;
    \item $R(e_{21}) =  -ve_{22}$, $R(ve_{11}) =  ve_{12}$, $R(ve_{21}) = e_{11}$;
    \item $R(e_{11}) =  e_{12}$, $R(e_{21}) = - e_{11}$, $R(ve_{21}) =  - ve_{11}$, $R(ve_{22}) =  - ve_{12}$;
    \item $R(e_{11}) =  e_{12}$, $R(e_{21}) = - e_{11}  -  ve_{11}$, $R(ve_{21}) =  - ve_{11}$, $R(ve_{22}) =  e_{12} - ve_{12}$;
    \item $R(e_{11}) = ve_{12}$, $R(e_{21}) =  -\alpha ve_{11}$, $R(ve_{21}) = e_{11}$, $R(ve_{22}) =  \alpha e_{12}$;
    \item $R(e_{11}) = -e_{12} + ve_{12}$, $R(e_{21}) =  e_{11} - \alpha ve_{11}$, $R(ve_{21}) = e_{11} + ve_{11}$, $R(ve_{22}) = \alpha e_{12} + ve_{12}$, $\alpha\neq -1$;
    \item $R(e_{22}) = ve_{12}$, $R(e_{21}) =  - \alpha ve_{11}$, $R(ve_{21}) =  - e_{11}$, $R(ve_{22}) =  \alpha e_{12}$.
\end{enumerate}
}

\textbf{Proof.} In \cite{Octonions} it was proved that there are only seven non-zero non-unital subalgebras $B$ in $\mathbb{O}$, up to action of automorphism. They are precisely the subalgebras from the statements of Lemmas~1--7. The theorem is proven.

\medskip\textbf{Corollary 6.} \textit{Let $R$ be the Rota-Baxter operator on the split Cayley-Dickson algebra $\mathbb{O}$. If a~field $F$ is quadratically closed, then, up to the action of automorphism, antiautomorphism and multiplication by a~scalar, $R$ acts on $\mathbb{O}$ in one of the following ways for some $\alpha\in F$ (unspecified basic elements $e_{ij}$, $ve_{ij}$ lie in $\ker (R)$):
\begin{enumerate}
    \item $R(e_{21}) = e_{12}$;
    \item $R(ve_{22}) = e_{12}$;
    \item $R(e_{21}) = e_{11}$;
    \item $R(e_{21})=e_{11}$, $R(e_{22})=e_{12}$;
    \item $R(e_{21})=-e_{11}$, $R(e_{11})=e_{12}$;
    \item $R(e_{21})=e_{11}$, $R(ve_{21})= e_{12}$;
    \item $R(ve_{11})=e_{11}$, $R(ve_{21})= e_{12}$;
    \item $R(ve_{11}) = e_{12}, R(ve_{21}) = e_{11}$;
    \item $R(ve_{21}) = e_{11}, R(ve_{22}) = e_{12}$; 
    \item $R(e_{21}) = ve_{12}$, $R(ve_{21}) = ve_{22}$;
    \item $R(e_{11}) = R(e_{12}) = -R(e_{21}) = - R(e_{22}) = ve_{22} + ve_{12}, R(ve_{11}) = - ve_{12}, R(ve_{21}) = ve_{12}$;
    \item $R(ve_{11}) = ve_{12}$, $R(ve_{21}) = ve_{22} + ve_{12}$;
    \item $R(ve_{11}) =  ve_{12}$, $R(ve_{21}) = \alpha ve_{22}$, $\alpha\neq 0$;
    \item $R(e_{21})= e_{12}$, $R(ve_{11})=ve_{12}$, $R(ve_{21})=ve_{22} + ve_{12}$;
    \item $R(e_{21}) = e_{12}$, $R(ve_{11})=ve_{12}$, $R(ve_{21})=ve_{22}$;
    \item $R(e_{21})= e_{12}$, $R(ve_{11})=ve_{12}$, $R(ve_{21})=ve_{22} + e_{12}$;
    \item $R(e_{21})=\alpha e_{12}$, $R(ve_{11})=ve_{12}$, $R(ve_{21})=ve_{22} + ve_{12} + e_{12}$, $\alpha\neq 0$;
    \item $R(e_{21}) =  ve_{12}$, $R(ve_{11}) = e_{11}$, $R(ve_{21}) = ve_{22}$;
    \item $R(e_{21}) =  e_{11}$, $R(ve_{11}) = \alpha ve_{12}$, $R(ve_{21}) = ve_{22}$, $\alpha\neq 0$;
    \item $R(e_{21}) =  -ve_{22}$, $R(ve_{11}) =  ve_{12}$, $R(ve_{21}) = e_{11}$;
    \item $R(e_{11}) =  e_{12}$, $R(e_{21}) = - e_{11}$, $R(ve_{21}) =  - ve_{11}$, $R(ve_{22}) =  - ve_{12}$;
    \item $R(e_{11}) =  e_{12}$, $R(e_{21}) = - e_{11}  -  ve_{11}$, $R(ve_{21}) =  - ve_{11}$, $R(ve_{22}) =  e_{12} - ve_{12}$;
    \item $R(e_{11}) = ve_{12}$, $R(e_{21}) =  - ve_{11}$, $R(ve_{21}) = e_{11}$, $R(ve_{22}) =  e_{12}$;
    \item $R(e_{11}) = -e_{12} + ve_{12}$, $R(e_{21}) =  e_{11} - \alpha ve_{11}$, $R(ve_{21}) = e_{11} + ve_{11}$, $R(ve_{22}) = \alpha e_{12} + ve_{12}$, $\alpha\neq -1$;
    \item $R(e_{22}) = ve_{12}$, $R(e_{21}) =  - ve_{11}$, $R(ve_{21}) =  - e_{11}$, $R(ve_{22}) =  e_{12}$.
\end{enumerate}
}

\medskip\textbf{Remark 1.} \textit{In the Corollary~6 operators (5), (21)--(24) are the ones, where $R^2\neq 0$, but $R^3 = 0$. Operators (1)--(4), (6)--(20) are the ones, where $R^2 = 0$.}

    \section{Acknowledgements}

The study was supported by a~grant from the Russian Science Foundation №~23-71-10005, https://rscf.ru/project/23-71-10005/

The author is very grateful to Vsevolod Gubarev for his advices in Rota-Baxter theory.

\end{document}